\documentclass[11pt,reqno]{amsart}
\usepackage{enumerate}
\usepackage{amsfonts}
\usepackage{amsmath}
\usepackage{amsthm}
\usepackage{amssymb}
\usepackage{latexsym}
\usepackage{multicol}
\usepackage{verbatim}
\usepackage{amscd}
\usepackage{hyperref}
\usepackage{pb-diagram}

\advance\textwidth by 1.2in \advance\oddsidemargin by -.6in
\advance\evensidemargin by -.6in
\newtheorem*{cor}{Corollary}
\newtheorem*{lem}{Lemma}
\newtheorem*{prop}{Proposition}
\allowdisplaybreaks[1]

\theoremstyle{definition}

\theoremstyle{definition}
\newtheorem{thm}{Theorem}

\newtheorem*{rem}{Remark}

\newenvironment{pf}{\proof}{\endproof}
\newcounter{cnt}

\makeatletter
\def\mydggeometry{\makeatletter\dg@YGRID=1\dg@XGRID=20\unitlength=0.003pt\makeatother}
\makeatother \theoremstyle{remark}

\numberwithin{equation}{section}

\DeclareMathOperator{\gldim}{gl.dim}

\DeclareMathOperator{\ad}{ad}
\let\bwdg\bigwedge
\def\bigwedge{{\textstyle\bwdg}}

\begin{document}

\newcommand{\thmref}[1]{Theorem~\ref{#1}}
\newcommand{\secref}[1]{Section~\ref{#1}}
\newcommand{\lemref}[1]{Lemma~\ref{#1}}
\newcommand{\propref}[1]{Proposition~\ref{#1}}
\newcommand{\corref}[1]{Corollary~\ref{#1}}
\newcommand{\remref}[1]{Remark~\ref{#1}}
\newcommand{\defref}[1]{Definition~\ref{#1}}
\newcommand{\er}[1]{(\ref{#1})}
\newcommand{\id}{\operatorname{id}}
\newcommand{\sgn}{\operatorname{sgn}}
\newcommand{\wt}{\operatorname{wt}}
\newcommand{\tensor}{\otimes}
\newcommand{\mof}{\operatorname{mod}_f}
\newcommand{\from}{\leftarrow}
\newcommand{\nc}{\newcommand}
\newcommand{\rnc}{\renewcommand}
\newcommand{\dist}{\operatorname{dist}}
\newcommand{\qbinom}[2]{\genfrac[]{0pt}0{#1}{#2}}
\nc{\cal}{\mathcal} \nc{\goth}{\mathfrak} \rnc{\bold}{\mathbf}
\renewcommand{\frak}{\mathfrak}
\newcommand{\supp}{\operatorname{supp}}
\renewcommand{\Bbb}{\mathbb}
\nc\bomega{{\mbox{\boldmath $\omega$}}} \nc\bpsi{{\mbox{\boldmath
$\Psi$}}}
 \nc\balpha{{\mbox{\boldmath $\alpha$}}}
 \nc\bpi{{\mbox{\boldmath $\pi$}}}

\newcommand{\lie}[1]{\mathfrak{#1}}
\makeatletter
\def\section{\def\@secnumfont{\mdseries}\@startsection{section}{1}%
  \z@{.7\linespacing\@plus\linespacing}{.5\linespacing}%
  {\normalfont\scshape\centering}}
\def\subsection{\def\@secnumfont{\bfseries}\@startsection{subsection}{2}%
  {\parindent}{.5\linespacing\@plus.7\linespacing}{-.5em}%
  {\normalfont\bfseries}}
\makeatother
\def\subl#1{\subsection{}\label{#1}}
 \nc{\Hom}{\operatorname{Hom}}
\nc{\End}{\operatorname{End}} \nc{\wh}[1]{\widehat{#1}}
\nc{\Ext}{\operatorname{Ext}} \nc{\ch}{\text{ch}}
\nc{\ev}{\operatorname{ev}} \nc{\Ob}{\operatorname{Ob}}
\nc{\soc}{\operatorname{soc}} \nc{\rad}{\operatorname{rad}}
\nc{\head}{\operatorname{head}}
\def\Im{\operatorname{Im}}
\def\gr{\operatorname{gr}}
\def\mult{\operatorname{mult}}

 \nc{\Cal}{\cal} \nc{\Xp}[1]{X^+(#1)} \nc{\Xm}[1]{X^-(#1)}
\nc{\on}{\operatorname} \nc{\Z}{{\bold Z}} \nc{\J}{{\cal J}}
\nc{\C}{{\bold C}} \nc{\Q}{{\bold Q}}
\renewcommand{\P}{{\cal P}}
\nc{\N}{{\Bbb N}} \nc\boa{\bold a} \nc\bob{\bold b} \nc\boc{\bold c}
\nc\bod{\bold d} \nc\boe{\bold e} \nc\bof{\bold f} \nc\bog{\bold g}
\nc\boh{\bold h} \nc\boi{\bold i} \nc\boj{\bold j} \nc\bok{\bold k}
\nc\bol{\bold l} \nc\bom{\bold m} \nc\bon{\bold n} \nc\boo{\bold o}
\nc\bop{\bold p} \nc\boq{\bold q} \nc\bor{\bold r} \nc\bos{\bold s}
\nc\bou{\bold u} \nc\bov{\bold v} \nc\bow{\bold w} \nc\boz{\bold z}
\nc\boy{\bold y} \nc\ba{\bold A} \nc\bb{\bold B} \nc\bc{\bold C}
\nc\bd{\bold D} \nc\be{\bold E} \nc\bg{\bold G} \nc\bh{\bold H}
\nc\bi{\bold I} \nc\bj{\bold J} \nc\bk{\bold K} \nc\bl{\bold L}
\nc\bm{\bold M} \nc\bn{\bold N} \nc\bo{\bold O} \nc\bp{\bold P}
\nc\bq{\bold Q} \nc\br{\bold R} \nc\bs{\bold S} \nc\bt{\bold T}
\nc\bu{\bold U} \nc\bv{\bold V} \nc\bw{\bold W} \nc\bz{\bold Z}
\nc\bx{\bold x}
\title[A family of Koszul algebras from simple Lie algebras]{A family of Koszul algebras arising from finite-dimensional representations of simple Lie algebras}
\author{Vyjayanthi Chari and  Jacob Greenstein}
\thanks{This work was partially supported by the NSF grants DMS-0500751 and DMS-0654421}
\address{Department of Mathematics, University of
California, Riverside, CA 92521.} \email{chari@math.ucr.edu}
 \email{jacob.greenstein@ucr.edu}\maketitle
\maketitle
\begin{abstract} Let $\lie g$ be a simple Lie algebra and let~$\bs^{\lie g}$ be the locally finite part of
the algebra of invariants $(\End_\bc\bv\otimes S(\lie g))^{\lie g}$ where $\bv$ is the direct sum of all simple finite-dimensional modules for $\lie g$ and $S(\lie g)$ is the symmetric algebra of $\lie g$. Given an integral weight $\xi$, let $\Psi=\Psi(\xi)$ be the subset of roots which have maximal scalar product with $\xi$. Given a dominant integral weight $\lambda$ and $\xi$ such that $\Psi$ is a subset of the positive roots we construct a finite-dimensional subalgebra $\bs^{\lie g}_\Psi(\le_\Psi\lambda)$ of $\bs^{\lie g}$ and prove that the algebra is Koszul of global dimension at most the cardinality of $\Psi$. Using this we then construct naturally an infinite-dimensional Koszul algebra of global dimension equal to the cardinality of $\Psi$. The results and the methods are motivated by the study of the category of finite-dimensional representations of the affine and quantum affine algebras.
\end{abstract}

\section*{Introduction}

This paper is motivated by the study of the category of finite-dimensional representations of both the classical and quantum loop algebras associated to a
simple Lie algebra $\lie g$. This category is not semisimple and thus it is natural to investigate its homological properties.
However, this category is both too large in the sense that it has too many simple objects and too small in the sense that it does not have enough projectives.
This means that one of the classical tools of representation theory, namely replacing a category by the category of modules over the endomorphism ring
of its projective generator, is not available to us. This tool  plays a very important role in the study of the category~$\cal O$ (cf.~\cite{BGG} and more
recently~\cite{BGS,BKM,FKM,Soe,Str} to name but a few)  and in many other situations as well. 
It is also  well-known that  the endomorphism ring of a projective generator
 is often a nice associative algebra. For instance, in the case of category~$\cal O$ that  algebra is  Koszul. 
Thus another motivation for this paper was to find an appropriate category of finite-dimensional modules which would allow us to use these methods.

There are many important families of irreducible representations of the quantum affine algebra, such as the Kirillov-Reshetikhin modules (\cite{Kir}) or more generally the minimal affinizations associated to dominant integral weights in the weight lattice of the simple Lie algebra (\cite{Ch1}), which on specializing to $q=1$  become
indecomposable modules for a truncated loop algebra (\cite{Ch2}). In most cases, they are in fact modules for $\lie g\tensor \bc[t^{\pm 1}]/((t-1)^2)$
and admit a natural $\bz_+$-grading given by powers of $(t-1)$.
In other words, some important families of modules can be regarded as $\bz_+$-graded modules for the finite-dimensional $\bz_+$-algebra $\lie g\ltimes\lie g_{\ad}$, where $\lie g_{\ad}$ has degree~1 and is an abelian ideal. It is therefore reasonable to restrict our attention to the subcategory $\cal G_2$ of graded modules
with finite-dimensional graded pieces over this quotient of the loop algebra.  This category still has infinitely many simple objects but they are now parametrized discretely and admit projective covers. In fact,
we show that every simple module has an explicit
projective resolution, coming from the Koszul complex of the symmetric algebra of~$\lie g$, and this  allows us to compute all~$\Ext$ spaces.

To construct interesting finite-dimensional associative algebras we have to pass to  Serre subcategories which contain only finitely many simple objects.  Using
the description of the
Kirillov-Reshetikhin modules in terms of generators and relations given in~\cite{Ch2, CMkir1}, one can naturally associate to such a module a  subset $\Psi_i$ of $R^+$. This subset  maximizes the linear
functional on the dual of a Cartan subalgebra of $\lie g$ given by the scalar product with the $i$th fundamental weight $\omega_i$. It is then natural to consider the case when
$\omega_i$ is replaced by an  arbitrary integral weight~$\xi$. We denote the corresponding set by~$\Psi(\xi)$.
These sets have many interesting combinatorial properties which are studied in~\cite{CDR}.
For instance, in some cases the maximal possible cardinality of
such a set equals the maximal dimension of a nilpotent abelian subalgebra of~$\lie g$ which were  computed already in~\cite{Mal}.

Suppose now that $\xi$ is such that $\Psi=\Psi(\xi)$ is contained in a fixed set of positive roots of~$\lie g$. This is always the case if~$\xi$
is dominant and in fact in this case $\Psi$ defines an abelian ideal in the corresponding Borel subalgebra. We define a partial order $\le_\Psi$ on the set of dominant integral weights which is a refinement of the usual partial order. Let $\cal G_2[\le_\Psi\lambda]$ be the subcategory of $\cal G_2$ consisting of objects whose irreducible constituents are in $\le_\Psi$. We show that $\cal G_2[\le_\Psi\lambda]$ has enough projectives
and is of global dimension at most equal to $|\Psi|$  and the bound is attained if $\lambda$ is sufficiently dominant.
The endomorphism ring of a projective generator of~$\cal G_2[\le_\Psi\lambda]$ admits a natural grading and we are able to prove by using the results of \cite{BGS} that this grading is Koszul.
We are also able to identify the endomorphism algebra as a subalgebra of $\lie g$-invariants in the tensor product algebra $\End_\bc\bv\otimes S(\lie g)$ where $\bv$ is the direct sum of all simple finite-dimensional modules for $\lie g$.
Finally we prove that these algebras ``approximate'' an  infinite-dimensional
algebra which is also Koszul and describe the  the corresponding Yoneda algebras explicitly. These algebras can be realized as path algebras of rather simple quivers with relations
(cf. example in~\ref{MR130} and~\cite{G}).

The paper is organized as follows. In Section~\ref{MR} we formulate the main results of the paper and briefly explain the strategy for proving them. In Section~\ref{CAT} we
study the fundamental properties of the category~$\cal G_2$ and a family of its Serre subcategories. In Section~\ref{HOM} we proceed to investigate the homological properties
of the category~$\cal G_2$. The next two sections are dedicated to studying the endomorphism algebras and their  quadratic duals.

\subsection*{Acknowledgments}
We are grateful to Arkady Berenstein, Bernhard Keller, Bertram Kostant, Bernard Leclerc and Sergey Loktev for many stimulating discussions.
The second author thanks Jonathan Brundan for some pointers to the appropriate results in~\cite{BGS}.

\section{Main results}\label{MR}

Throughout this paper $\bz_+$ stands for the set of non-negative integers and~$\bc$ denotes the field of complex numbers. All algebras and
vector spaces, as well as $\Hom$ and tensor products, are considered over~$\bc$ unless specified otherwise. If~$A$ is an associative algebra,
we denote~$A^{op}$ its opposite algebra.

\subsection{}\label{MR10} Let~$\lie g$ be a complex finite dimensional simple Lie algebra. Fix a Cartan subalgebra~$\lie h$ and let~$R\subset P\subset \lie h^*$ be, respectively,
the set of roots and the weight lattice of~$\lie
g$ with respect to~$\lie h$. Let $(\cdot,\cdot)$ be the
non-degenerate symmetric bilinear form on $\lie h^*$ induced by
the restriction of the Killing form of $\lie g$ to $\lie h$. Set
$I=\{1,\dots ,\dim\lie h\}$ and let $\alpha_i$, $i\in I$ be a set
of simple roots and $\omega_i$, $i\in I$ the corresponding
fundamental weights. Let $R^+$ be the corresponding set of positive
roots and let $P^+$ be the $\bz_+$-span of the fundamental weights.

\subsection{}\label{MR30} For~$\lambda\in P^+$, let~$V(\lambda)$ be the unique, up to  isomorphism,
simple finite dimensional  $\lie g$-module with  highest weight~$\lambda$ and let $V(\lambda)^*=\Hom(V(\lambda),\bc)$ with
its standard $\lie g$-module structure. Set
$$\bv=\bigoplus_{\lambda\in P^+} V(\lambda),\qquad \bv^{\circledast}=\bigoplus_{\lambda\in P^+} V(\lambda)^*.$$
Note that~$\bv^{\circledast}\cong \bv$ as a $\lie g$-module.
Define an associative algebra structure on $\bv^{\circledast}\otimes\bv$ by extending linearly
$$(f\tensor v)(g\otimes w)=g(v)f\otimes w,\qquad v,w\in \bv,\, f,g\in \bv^{\circledast}.
$$ 
The natural embedding
$\bv^{\circledast}\otimes \bv\to\End\bv$ (respectively, $\bv^{\circledast}\otimes \bv\to\End\bv^\circledast$) of
$\lie g$-modules given by extending $f\tensor v\mapsto (w\mapsto f(w)v)$ (respectively, $f\tensor v\mapsto (g\mapsto g(v)f$)) 
is then an anti-homomorphism of associative algebras. In particular, for all~$\lambda\in P^+$ we have isomorphisms of associative algebras
$$
V(\lambda)^*\tensor V(\lambda)\to (\End V(\lambda))^{op},\qquad V(\lambda)^*\tensor V(\lambda)\to (\End V(\lambda)^*)^{op}.
$$
The preimage of the identity element in~$\End V(\lambda)^{op}$ (or~$(\End V(\lambda)^*)^{op}$) is the canonical $\lie g$-invariant 
element~$1_\lambda$ of~$V(\lambda)^*\tensor V(\lambda)$.
\begin{lem}
For
$\lambda,\mu\in P^+$ we have
$$1_\lambda1_\lambda=1_\lambda,\qquad 1_\lambda1_\mu=0,\quad \mu\ne\lambda$$ and
$$1_\lambda(\bv^{\circledast}\otimes \bv)=V(\lambda)^*\otimes \bv,\qquad (\bv^\circledast\tensor \bv) 1_\mu=\bv^{\circledast}\tensor V(\mu).$$
\end{lem}
\begin{pf}
Observe first that~$(f\tensor v)(g\tensor w)=0$ if~$v\in V(\nu)$, $g\in V(\xi)^*$ with~$\xi\not=\nu$. 
Write~$1_\lambda=\sum_i \xi_i\tensor u_i$, $\xi_i\in V(\lambda)^*$, $u_i\in V(\lambda)$.
Then we have $\sum_i \xi_i(v) u_i=v$ for all~$u\in V(\lambda)$ and $\sum_i \xi_i f(u_i)=f$ for all~$f \in V(\lambda)^*$.
It follows that~$1_\lambda(f\tensor v)=f\tensor v$ for all~$f\in V(\lambda)^*$, $v\in\bv$ and that~$(g\tensor u)1_\lambda=
g\tensor u$ for all~$g\in \bv^\circledast$, $u\in V(\lambda)$. Since~$\bv^\circledast\tensor \bv=\bigoplus_{\nu,\xi\in P^+} V(\nu)^*\tensor V(\xi)$,
the assertions follow.
\end{pf}

\subsection{}\label{MR50} Let $A$ be any associative algebra with unity $1_A$ and consider $ \ba = A\otimes (\bv^{\circledast}\otimes\bv)$.
This is obviously an associative algebra. If  $A$ is a $\lie
g$-module algebra, that is the multiplication map $A\tensor A\to A$ is a homomorphism of $\lie g$-modules, then the usual tensor product action defines a
$\lie g$-module structure on $\ba$ and the
algebra multiplication is again a morphism of $\lie g$-modules. By
abuse of notation we let $1_\lambda$ also denote the
idempotent $1_A\otimes 1_\lambda$ in~$\ba$. For $\lambda,\mu\in P^+$ we
have 
$$
1_\lambda \ba 1_\mu= A\otimes V(\lambda)^*\otimes
V(\mu).
$$ 
If the algebra $A$ has  a $\bz_+$-grading $A=\bigoplus_{k\in\bz_+} A[k]$ which is
compatible with the $\lie g$-action, then $\ba$ has a natural
grading given by
\begin{gather*} \ba[k]=A[k]\otimes (\bv^{\circledast}\otimes\bv).
\end{gather*}
\noindent{\it From now on we shall assume  that $A$ is a
$\bz_+$-graded associative $\lie g$-module algebra with unity and
also that $A[0]=\bc 1_A$, $\dim A[k]<\infty$, $k\in\bz_+$.}

\subsection{}\label{MR70} Given~$\Psi\subset R^+$, let~$\le_\Psi$ be the
partial order on~$P^+$ given by:  $$\lambda\le_\Psi \mu\ \iff
\mu-\lambda\in \bz_+ \Psi,$$  where $\bz_+\Psi$ is the
non-negative integer linear span of the elements of $\Psi$. This
is a refinement of the usual partial order~$\le=\le_{R^+}$ on $P^+$. Given $\lambda,\mu\in P^+$, set
\begin{align*}
&\le_\Psi\lambda= \{ \nu\in P^+\,:\, \nu\le_\Psi\lambda\},
\qquad
\lambda\le_\Psi=\{ \nu\in P^+\,:\, \lambda\le_\Psi\nu\},\\
&[\mu, \lambda]_\Psi = (\le_\Psi\lambda)\cap (\mu\le_\Psi).
\end{align*}
The first set and hence the third are finite subsets of $P^+$.

For all $\lambda,\mu\in P^+$ with
$\lambda\le_\Psi\mu$, define
$$d_\Psi(\lambda,\mu)=\min\{\sum_{\beta\in\Psi} m_\beta:
\mu-\lambda=\sum_{\beta\in \Psi} m_\beta\beta,\
m_\beta\in\bz_+\}.$$
Given $\lambda,\mu\in P^+$ with $\lambda\le_\Psi\mu$, set
$$\ba_\Psi(\lambda,\mu)=1_\lambda\ba[d_\Psi(\lambda,\mu)]1_\mu$$
and for $F\subset P^+$, set
$$\ba_\Psi(F)=\bigoplus_{\lambda,\mu\in F\,:\,\lambda\le_\Psi \mu}
\ba_\Psi(\lambda,\mu).$$  Clearly $\ba_\Psi(\lambda,\mu)$ and $\ba_\Psi(F)$ are $\lie g$-submodules of $\ba$. The following Lemma is easily checked.
\begin{lem} Suppose that $\Psi\subset R^+$ is such that for all
$\lambda,\mu,\nu\in P^+$ with $\lambda\le_\Psi\mu\le_\Psi\nu$, we
have $$d_\Psi(\lambda,\mu)+d_\Psi(\mu,\nu)=d_\Psi(\lambda,\nu).$$
Then for all $F\subset P^+$ the $\lie g$-submodule $\ba_\Psi(F)$
is a graded subalgebra of $\ba$.\hfill\qedsymbol
\end{lem}
It is easy to construct examples of sets $\Psi$ satisfying the conditions of the Lemma. Thus for $ \xi\in P$, let $$\max \xi=\max\{(\alpha,\xi): \alpha\in R\} ,\qquad\Psi(\xi)=\{\alpha\in R: (\xi,\alpha)=\max\xi\}.$$
Clearly  if $\xi\ne 0$, we have $$\max\xi>0,\qquad (\Psi(\xi)+\Psi(\xi))\cap (R\cup\{0\})=\emptyset.$$

Let~$\Psi=\Psi(\xi)$ and assume that~$\Psi\subset R^+$ (note that this holds if~$\xi\in P^+$). Then we see that for  all $\lambda,\mu,\nu\in P^+$ with $\lambda\le_\Psi\mu\le_\Psi\nu$ we have $$d_\Psi(\lambda,\mu)+d_\Psi(\mu,\nu)=d_\Psi(\lambda,\nu).$$
\subsection{}\label{MR90}
 Let $\ba^{\lie g}$ be the submodule of $\lie g$-invariants of $\ba$. For $\lambda,\mu\in P^+$, $F\subset P^+$, and $\Psi\subset R^+$ set $$\ba^{\lie g}_\Psi(\lambda,\mu) =(\ba_\Psi(\lambda,\mu))^{\lie g},\qquad\ba_\Psi^{\lie g}(F)=(\ba_\Psi(F))^{\lie g}.$$ The following is clear.
 \begin{prop}
 \begin{enumerate}[{\rm(i)}]
\item\label{MR90.i} The submodule  $\ba^{\lie g}$ is a graded subalgebra of $\ba$ and $$\ba^{\lie g}[k]=(\ba[k])^{\lie g}.$$
 \item\label{MR90.ii} If $\Psi$ is a subset of $R^+$ such that for all
$\lambda,\mu,\nu\in P^+$ with $\lambda\le_\Psi\mu\le_\Psi\nu$, we
have $d_\Psi(\lambda,\mu)+d_\Psi(\mu,\nu)=d_\Psi(\lambda,\nu)$, then $\ba^{\lie g}_\Psi(F)$ is a graded subalgebra of $\ba^{\lie g}$ for all $F\subset P^+$.
\end{enumerate}
 \end{prop}

\subsection{}
We denote $\bt$ (respectively, $\bs$, $\be$) the algebra~$\ba$ with~$A=T(\lie g)$ (respectively, $A=S(\lie g)$, $A=\bigwedge \lie g$),
the $\lie g$-module structure on~$A$ being given by the usual diagonal action.
 Our first result is the following.
\begin{thm}\label{kosalg} Let $\xi\in P$ be such that $\Psi=\Psi(\xi)\subset R^+$.  Let $A$ be  either  $S(\lie g)$ or $\bigwedge \lie g$.

\begin{enumerate}[{\rm(i)}]
\item\label{kosalg.i} Let $\mu,\nu\in P^+$. The  subalgebras $\ba_\Psi^{\lie g}(\le_\Psi\nu)$,
$\ba_\Psi^{\lie g}(\mu\le_\Psi)$ and $\ba_\Psi^{\lie g} ([\mu,\nu]_\Psi)$ of
$\ba_\Psi^{\lie g}$ are Koszul.
For $A=S(\lie g)$ all these subalgebras have global dimension at most~$|\Psi|$ and
the   bound is attained for some $\mu',\nu'\in P^+$ with $\mu'\le_\Psi\nu'$.

\item\label{kosalg.ii} The algebra $\ba_\Psi^{\lie g}(P^+)$ is Koszul.
Moreover, $\bs_\Psi^{\lie g}$ has left global dimension $|\Psi|$ and its quadratic dual is $(\be_\Psi^\lie g)^{op}$.
\end{enumerate}
\end{thm}

\begin{rem} An alternative characterization of the sets $\Psi(\xi)$ can be found in \cite{CDR} together with a complete list of such sets for~$\lie g$ of classical types. In particular, one sees that as the rank of the Lie algebra varies one can find sets $\Psi(\xi)$ of arbitrary cardinality. As a consequence, we see that the Theorem implies the existence of an infinite-dimensional Koszul algebra of global dimension $k$ for any $k\ge 1$.
\end{rem}

\subsection{}\label{MR130}
As an example of our construction, we produce an algebra of left global dimension~$2$. Assume that $\lie g$ is of rank greater than two and that $\lie g$ not of type $A$ or $C$. Let $\theta\in R^+$ be the highest root and choose $i_0\in I$ such that $(\theta,\alpha_{i_0})>0$. Then $2\theta-\alpha_{i_0}\in P^+$ and $\Psi(2\theta-\alpha_{i_0})=\{\theta,\theta-\alpha_{i_0}\}.$
In this case it can be shown that~$\bs_\Psi^{\lie g}$ is isomorphic to an infinite  direct sum of copies of the algebra $\mathfrak B$ defined as follows.
Consider  the following translation quiver
$$
\begin{CD}
\scriptstyle
(0,0) @>>> \scriptstyle(0,1) @>>> \scriptstyle(0,2) @>>> \scriptstyle(0,3) @>>>\cdots\\
@. @VVV @VVV @VVV\\
@. \scriptstyle(1,0) @>>> \scriptstyle(1,1) @>>>\scriptstyle(1,2) @>>> \cdots\\
@. @. @VVV @VVV\\
@. @. \scriptstyle(2,0) @>>> \scriptstyle(2,1) @>>> \cdots\\
@. @. @. @VVV \\
@. @. @. \scriptstyle(3,0) @>>> \cdots
\end{CD}
$$
the translation map being~$\tau(m,n)=(m-1,n)$, $m,n\in\bz_+$, $m>0$. Then $\mathfrak B$ is the quotient of the path algebra of the quiver  by the mesh relations. Thus, this algebra is an  infinite dimensional analogue
of the Auslander algebra of the path algebra of type~$\mathbb A_\infty$.
Another example of a finite dimensional algebra from this family was already constructed in~\cite{CG}. Other examples will
 appear in \cite{G}.

 \subsection{}\label{MR150} In
the rest of the section we explain the motivation for the
definition of the algebras $\ba^{\lie g}_\Psi$ and our strategy
for proving Theorem \ref{kosalg}.  We  prove that $\bs^{\lie
g}_\Psi(\le_\Psi\lambda)$ is isomorphic to the endomorphism
algebra of a projective generator of a category of modules for a
certain finite-dimensional algebra $\lie g\ltimes\lie g_{\ad}$
defined as follows. As a vector space $$\lie g\ltimes\lie
g_{\ad}=\lie g\oplus\lie g,$$ and the Lie bracket is given by,
$$[(x,y),(x',y')]=([x,x'],[x,y']+[y,x']).$$ In particular if we
identify $\lie g$ (respectively, $\lie g_{\ad}$) with the subspace
$\{(x,0):x\in\lie g\}$ (respectively, $\{(0,y): y\in\lie g\}$), then $\lie
g_{\ad}$ is an abelian Lie ideal in $\lie g\ltimes\lie g_{\ad}$.
Define a $\bz_+$-grading on $\lie g\ltimes \lie g_{\ad}$ by
requiring the elements of $\lie g$ to have degree zero and elements
of $\lie g_{\ad}$ to have degree one. Then the universal enveloping
algebra  $\bu(\lie g\ltimes\lie g_{\ad})$ is a $\bz_+$-graded
algebra and as a trivial consequence of the PBW theorem, there is
an isomorphism of vector spaces $$\bu(\lie g\ltimes\lie
g_{\ad})\cong S(\lie g)\tensor \bu(\lie g).$$
Our main motivation for the study of $\lie g\ltimes\lie g_{\ad}$ as a $\bz_+$-graded Lie algebra stems from the easy observation that
$$
\lie g\ltimes\lie g_{\ad}\cong \lie g\tensor \bc[t,t^{-1}]/( (t-1)^2)
$$
as $\bz_+$-graded Lie algebras, the right hand side being graded by powers of~$(t-1)$.

\subsection{} Let $\cal G_2$ be the category whose objects are
 $\bz_+$-graded $\lie g\ltimes\lie g_{\ad}$-modules with finite-dimensional graded pieces and where the
 morphisms are $\lie g\ltimes\lie g_{\ad}$-module maps which preserve the grading.
 In other words a $\lie g\ltimes\lie g_{\ad}$-module $V$ is an object of $\cal G_2$
 if and only if 
\begin{gather*} 
V=\bigoplus_{k\in\bz_+} V[k],\qquad\dim V[k]<\infty,\\
 \lie g V[k]\subset V[k],\qquad \lie g_{\ad }V[k]\subset
 V[k+1].
\end{gather*}
If $V,W\in\Ob\cal G_2$, then
  $$\Hom_{\cal G_2}(V,W)=\{f\in\Hom_{\lie g\ltimes\lie g_{\ad}}(V,W): f(V[k])\subset W[k]\}.$$

 We prove that simple objects in this category are parametrized by
  the set $\Lambda=P^+\times \bz_+$. For $(\lambda,r)\in \Lambda$ let $V(\lambda,r)$
  be an element in the corresponding isomorphism class. We then prove
\begin{prop}\label{extsimple}
For all~$j\ge 0$, $(\mu,s),(\lambda,r)\in\Lambda$,
$$
\Ext^j_{{\cal G}_2}(V(\mu,s),V(\lambda,r))\cong \begin{cases}
\Hom_{\lie g}(\textstyle\bigwedge^{j} \lie g\tensor
V(\mu),V(\lambda)),& j=r-s, \\ 0,& {\text{otherwise}.}
\end{cases}
$$
\end{prop}
\subsection{}\label{MR170} Suppose now that we have a subset $\Psi$ of $R^+$ as in Theorem~\ref{kosalg}.
 Given $\lambda\in P^+$, let $\cal G_2[\le_\Psi\lambda]$ be the full subcategory of $\cal G_2$
  consisting of  finite-dimensional objects $V\in\cal G_2$  satisfying the following: if $V(\nu,s)$
   is a simple constituent of $V$, then $\nu\le_\Psi\lambda$ and $s=d_\Psi(\nu,\lambda)$.
   In particular, $\cal G_2[\le_\Psi\lambda]$ has only finitely many simple objects and
    is a Serre subcategory of~$\cal G_2$.
\begin{thm}\label{equivalence} The category $\cal G_2[\le_\Psi\lambda]$ has enough projectives.
 If $V(\nu_k,s_k)$, $k=1,2$ are simple objects in $\cal G_2[\le_\Psi\lambda]$ then
$$
\Ext^j_{{\cal G}_2}(V(\nu_1,s_1),V(\nu_2,s_2))\cong\Ext^j_{{\cal G}_2[\le_\Psi\lambda]}(V(\nu_1,s_1),V(\nu_2,s_2)).$$
If $P(\le_\Psi\lambda)$ is a projective generator of $\cal
G_2[\le_\Psi\lambda]$ we have an isomorphism of $\bz_+$-graded
associative algebras
$$
\End_{\cal G_2[\le_\Psi\lambda]} P(\le_\Psi\lambda)\cong \bs^{\lie
g}_\Psi(\le_\Psi \lambda)^{op}
$$
and the category~$\cal G_2[\le_\Psi\lambda]$ is equivalent to the
category of left finite-dimensional $\bs^{\lie g}_\Psi(\le_\Psi
\lambda)$-modules. In particular, $\bs^{\lie g}_\Psi(\le_\Psi\lambda)$ is
Koszul and $\gldim \bs^{\lie g}_\Psi(\le_\Psi\lambda)\le |\Psi|$.
\end{thm}
\subsection{} Once~\thmref{equivalence} is established, we then prove that $\bs^{\lie g}_\Psi$ is
 quadratic and that its Koszul complex is exact. Further we prove that its quadratic dual
  is $\be^{\lie g}_\Psi{}^{op}$ which implies that $\be^{\lie g}_\Psi$ is also Koszul.
   In particular, this also proves that $\be^{\lie g}_\Psi(\le_\Psi\lambda)$
is isomorphic to the Yoneda algebra of~$\bs^{\lie
g}_\Psi(\le_\Psi\lambda)$ which then establishes~\thmref{kosalg}.

{\em In the rest of the paper, we will identify the algebra $\ba$
with $\bv^{\circledast}\tensor A\tensor \bv$ under the natural isomorphism of~$\lie g$-modules. The algebra structure induced on $\bv^{\circledast}\tensor A\tensor \bv$ is given by $$(f\otimes a\otimes v)(g\otimes b\otimes w)= g(v)f\otimes ab\otimes w,\qquad a,b\in A,\, f,g\in\bv^{\circledast},\, v,w\in\bv .$$}

\section{The categories \texorpdfstring{$\cal G_2$ and $\cal G_2[\Gamma]$}{{\em G}\_2}}\label{CAT}
\subsection{}\label{PRE10}
Given $\alpha\in
R$ denote by $\lie g_\alpha$ the corresponding root space. The
subspaces $\lie n^\pm=\bigoplus_{\alpha\in R^+}\lie g_{\pm\alpha}$
are Lie subalgebras of $\lie g$. Fix a Chevalley basis
$x^\pm_\alpha$, $\alpha\in R^+$, $h_i$, $i\in I$ of $\lie g$ and
for $\alpha\in R^+$, set~$h_\alpha=[x_\alpha,x_{-\alpha}]$.

Let $\cal F(\lie g)$ be the category of finite-dimensional $\lie
g$-modules with the morphisms being maps of  $\lie g$-modules.  We
write $\Hom_{\lie g}$ for~$\Hom_{\cal F(\lie g)}$. If $V\in\Ob\cal
F(\lie g)$, denote
$$
V^{\lie g}=\{v\in V: \lie g v=0\}
$$
the subspace of~$\lie g$-invariants in~$V$.
Recall that
$$V=\bigoplus_{\lambda\in\lie h^*} V_\lambda,\qquad V_\lambda=\{v\in
V:hv=\lambda(h)v,\,\forall\, h\in\lie h\}.$$
Set
$$V^+=\{v\in
V\,:\, \lie n^+v=0\},\qquad V^+_\lambda=V^+\cap V_\lambda.$$
The category~$\cal F(\lie g)$ is semi-simple, i.e.
any object in~$\cal F(\lie g)$ is isomorphic to a direct sum of
simple modules.

Given $\lambda\in P^+$ let $v_\lambda\in V(\lambda)$ be such that
$$
\lie n^+ v_\lambda=0,\quad hv_\lambda=\lambda(h)v_\lambda,\quad
(x^-_{\alpha_i})^{\lambda(h_i)+1}v_\lambda =0,
$$
for all~$h\in\lie h$, $i\in I$.

We shall use the following standard
results in the course of the paper (cf.~\cite{PRV} for \eqref{PRE10.iv}).

\begin{lem} Let $\lambda,\mu\in P^+$, $V\in\Ob\cal F(\lie g)$. Then
\begin{enumerate}[{\rm(i)}]
\item\label{PRE10.iii} $\dim\Hom_{\lie
g}(V(\lambda), V)=\dim V^+_\lambda$,
\item\label{PRE10.iv} As vector spaces, we have
\begin{equation*}
\Hom_{\lie g}(V(\lambda), V\tensor V(\mu))\cong \{ v\in V_{\lambda-\mu}\,:\, (x_{\alpha_i}^+)^{\mu(h_i)+1}v=0 = (x_{\alpha_i}^-)^{\lambda(h_i)+1}v\}.\tag*{\qedsymbol}
\end{equation*}
\item\label{PRE10.v} Let $U,V,W\in\Ob\cal F(\lie g)$. The canonical  map $U^*\otimes V\to \Hom_\bc(U,V)$ is an isomorphism of $\lie g$-modules and we have $(U^*\otimes V)^{\lie g}\cong \Hom_{\lie g}(U,V)$.
The natural $\lie g$-module map $\Hom_\bc(U,V)\otimes \Hom_\bc(V,W)\to \Hom_\bc (U,W)$ given by composition and its restriction to $\Hom_{\lie g}(U,V)\otimes \Hom_{\lie g}(V,W)\to \Hom_{\lie g} (U,W)$ induces the  $\lie g$-module map $(U^*\otimes V)\otimes (V^*\otimes W)\to (U^*\otimes W)$ given by $$(f\otimes v)\otimes (v^*\otimes w)\mapsto v^*(v)f\otimes w,$$ and the restriction $(U^*\otimes V)^{\lie g}\otimes (V^*\otimes W)^{\lie g}\to (U^*\otimes W)^{\lie g}$
\end{enumerate}
\end{lem}

\subsection{}\label{CAT40}
Set $\Lambda=P^+\times\bz_+$. Given $(\lambda,r)\in\Lambda$, define
$V(\lambda,r)\in\Ob\cal G_2$ by $$V(\lambda,r)[s]=0,\qquad
r\ne s,\qquad V(\lambda,r)[r]\cong_{\lie g} V(\lambda),$$ with $\lie
g_{\ad}V(\lambda,r)=0$.   Observe that $V(0,0)\cong\lie\bc$ is the
trivial $\lie g\ltimes\lie g_{\ad}$-module. Set
\begin{equation}\label{CAT40.20}
 P(\lambda,r)=\bu(\lie
g\ltimes \lie g_{\ad})\otimes_{\bu(\lie g)} V(\lambda,r).
\end{equation}
Then $P(\lambda,r)\in\Ob\cal G_2$. It is immediate from the PBW
theorem that we have an isomorphism of $\bz_+$-graded $\lie
g\ltimes\lie g_{\ad}$-modules
$$
P(\lambda,r)\cong S(\lie g_{\ad})\otimes V(\lambda,r),\qquad
P(\lambda,r)[k]\cong_{\lie g} S^{k-r}(\lie g_{\ad})\otimes V(\lambda,r).
$$
For $V\in\Ob\cal G_2$ and $(\lambda,r)\in\Lambda$, set
$$[V:V(\lambda,r)]=\dim \Hom_{\lie g}(V(\lambda), V[r]).$$ If $V$
is finite-dimensional, then $[V:V(\lambda,r)]$ is just the
multiplicity of $V(\lambda,r)$ in  a Jordan-Holder series for $V$.
We have
\begin{equation}\label{CAT40.30}
[P(\lambda,r):V(\mu,s)]=\dim\Hom_{\lie g}(V(\mu),S^{s-r}(\lie g_{\ad})\otimes V(\lambda)).
\end{equation}

The next proposition is easily established along the lines
of~\cite[Proposition~2.1]{CG}. We include a short sketch of the
proof for the reader's convenience.
\begin{prop} Let $(\lambda,r)\in\Lambda$.
\begin{enumerate}[{\rm(i)}]
\item\label{CAT40.i}  The object $V(\lambda,r)$ is the unique
simple quotient of  $P(\lambda,r)$ and hence $P(\lambda,r)$ is its
projective cover in  $\cal G_2$. Moreover, the kernel of the
canonical morphism $P(\lambda,r)\to V(\lambda,r)$ is generated by
$P(\lambda,r)[r+1]$.

\item\label{CAT40.ii} The isomorphism classes of
simple objects in $\cal G_2$ are parametrized by pairs
$(\lambda,r)\in\Lambda$ and we have
\begin{align*}
&\Hom_{\cal G_2}(V(\lambda,r), V(\mu,s))=0,\qquad (\lambda,r)\ne (\mu,s),\\
&\Hom_{\cal G_2}(V(\lambda,r),V(\lambda,r))\cong\bc.\end{align*}

\item\label{CAT40.i5}  $P(\lambda,r)$ is the $\lie g\ltimes\lie
g_{\ad}$-module generated by an element $v_\lambda$ with
relations
$$(\lie n^+)v_\lambda =0,\qquad
hv_\lambda=\lambda(h)v_\lambda, \qquad
(x^-_{\alpha_i})^{\lambda(h_i)+1}v_\lambda=0,$$ for all $h\in\lie
h$ and $i\in I$. 

\item\label{CAT40.iv}If  $V$ is  concentrated in degree
$k$  for some $k\in\bz_+$, then $V$ is semi-simple and $P(V)=
S(\lie g_{\ad})\otimes V\in\Ob\cal G_2$  is the projective cover
of $V$ in~$\cal G_2$. Moreover, if $W$ is also concentrated in degree~$k$, then $P(V\oplus W)\cong P(V)\oplus P(W)$.
 \end{enumerate}
\end{prop}
\begin{pf} It is clear that $V(\lambda,r)$ is simple. To prove that $P(\lambda,r)$ is the
projective cover of $V(\lambda,r)$, note that it is indecomposable
since it is generated by $P(\lambda,r)[r]$ which is isomorphic to
the simple $\lie g$-module $V(\lambda)$. The fact that it is
projective is standard. To prove~\eqref{CAT40.ii}, it suffices to note that
any simple $\cal G_2$-module $V$ must satisfy $V[k]\ne 0$ for at
most one $k=r\in\bz_+$. In that case, $\lie g_{\ad}( V[r])\subset
V[r+1]=0$ and hence $V[r]\cong V(\lambda)$ for some $\lambda\in
P^+$.  The other parts of the proposition are now straightforward
and we omit the details.
\end{pf}

\subsection{}\label{CAT500}
{\em From now on we fix $\xi\in P\setminus\{0\}$ and assume that
$\Psi= \Psi(\xi)\subset R^+$ where $\Psi(\xi)=\{\alpha\in R:
(\alpha,\xi)=\max\xi\}$.} We will use the following property of
these sets repeatedly in the course of the paper.
\begin{lem}
Let~$\Psi=\Psi(\xi)$. Suppose that
$$
\sum_{\alpha\in R} m_\alpha \alpha=\sum_{\beta\in\Psi} n_\beta \beta,\qquad m_\alpha,n_\beta\in\bz_+.
$$
Then
\begin{equation}\label{CAT500.10}
\sum_{\beta\in\Psi} n_\beta\le \sum_{\alpha\in R} m_\alpha
 \end{equation}
with equality if and only if $m_\alpha=0$ for all~$\alpha\notin\Psi$. 
\end{lem}
\begin{pf}
 We have
\begin{equation*}
\max\xi \sum_{\beta\in\Psi} n_\beta=\sum_{\beta\in\Psi} n_\beta
(\beta,\xi)= \sum_{\alpha\in R} m_\alpha (\alpha,\xi)\le \max \xi
\sum_{\alpha\in R} m_\alpha.
\end{equation*}
The inequality in \eqref{CAT500.10} follows since~$\max\xi>0$, while the   equality holds if and only if
$$
\sum_{\alpha\in R} m_\alpha( \max\xi-(\alpha,\xi))=0.
$$
Since~$m_\alpha\in\bz_+$ this implies that equality holds if and only if~$m_\alpha=0$
unless~$(\alpha,\xi)=\max\xi$.
\end{pf}
\begin{rem} In \cite{CDR} it is shown that the converse of the
Lemma is true as well and that in fact if we let $\rho_\xi$ be the
sum of elements in $\Psi(\xi)$, then $\Psi(\xi)=\Psi(\rho_\xi)$.
\end{rem}

\subsection{}\label{CAT530} Let  $\Gamma$ be  a subset of $\Lambda$ and $\cal
G_2[\Gamma]$ be  the full subcategory of~$\cal G_2$  consisting of
objects $V$ such that
$$
[V: V(\mu,s)]\ne 0\, \implies \, (\mu,s)\in\Gamma.
$$
Let $V\to V_\Gamma$ be the functor from $\cal G_2\to\cal
G_2[\Gamma]$ defined by requiring $V_\Gamma$ to be the maximal
$\cal G_2[\Gamma]$-subobject of $V$.  In general one cannot
say much about these functors. However, in the next proposition
we shall see that for some special choices of $\Gamma$ and $V$
the module $V^\Gamma:= V/V_{\Lambda\setminus\Gamma}$ is an object in~$\cal G_2[\Gamma]$.

Given $\mu\le_\Psi\lambda\in
P^+$, let
\begin{align*}
&\Lambda(\le_\Psi\lambda)=\{(\nu, d_\Psi(\nu,\lambda)): \nu\in P^+,\,\nu\le_\Psi\lambda\},\\
&\Lambda([\mu,\lambda]_\Psi)=\{(\nu,d_\Psi(\nu,\lambda))\,:\,\nu\in P^+,\,\mu\le_\Psi\nu\le_\Psi\lambda\}.
\end{align*}
The set $\Lambda(\le_\Psi\lambda)$ (and hence~$\Lambda([\mu,\lambda]_\Psi)$) is finite and we denote the
corresponding category by $\cal G_2[\le_\Psi\lambda]$ (respectively, by $\cal G_2[[\mu,\lambda]_\Psi]$).

\begin{lem}
Let $\lambda'\le_\Psi \lambda\in P^+$. Set
$\Gamma=\Lambda(\le_\Psi\lambda)$ or~$\Gamma=\Lambda([\lambda',\lambda]_\Psi)$ and let $(\mu,r)\in\Gamma$.
We have
$$
[P(\mu,r)_{\Lambda\setminus\Gamma}:V(\nu,s)]=
\begin{cases}
[P(\mu,r):V(\nu,s)],& (\nu,s)\notin\Gamma,\\
0,&(\nu,s)\in\Gamma
\end{cases}
$$
In particular, $P(\mu,r)^\Gamma$ is an object in~$\cal G_2[\Gamma]$
and $[P(\mu,r)^\Gamma:V(\nu,s)]=[P(\mu,r):V(\nu,s)]$ for all~$(\nu,s)\in\Gamma$.
\end{lem}
\begin{pf} Set $P=P(\mu,r)$ and  for any subset $\Gamma'$ of $\Lambda$, let $P(\Gamma')$
be the $\lie g$-submodule of $P$ generated by elements of the subspace
$$\{v\in P[s]_\nu^+:
(\nu,s)\in\Gamma'\}.$$ The Lemma obviously follows if we
prove that $P(\Lambda\setminus\Gamma)$ is an object in $\cal G_2$,
i.e. that it is a $\bz_+$-graded  $\lie g\ltimes\lie
g_{\ad}$-module in which case we have
$P(\mu,r)_{\Lambda\setminus\Gamma}=P(\Lambda\setminus\Gamma)$. Let
$v_\nu\in P[s]^+_\nu$ for some $(\nu,s)\in\Lambda\setminus \Gamma$
and let $V=\bu(\lie g)v_\nu$. Consider the map $\lie
g_{\ad}\otimes V\to P[s+1]$ given by $x_{\ad}\otimes v\to x_{\ad}
v$ and suppose that its image  has a non-zero projection onto some $\lie
g$-module $\bu(\lie g)v_{\zeta}$, where $v_\zeta\in
P[s+1]_\zeta^+$. It suffices to prove that   $(\zeta,s+1)\in
\Lambda\setminus \Gamma$. Note  that
 either $\zeta=\nu$ or $\zeta=\nu-\beta_0$ for some
$\beta_0\in R$.  Assume for a contradiction that $(\zeta,s+1)\in\Gamma$.

Since $(\mu,r), (\zeta,s+1)\in\Gamma$, we can write $$\lambda-\mu=\sum_{\alpha\in\Psi} n_\alpha\alpha,\qquad \lambda-\zeta=\sum_{\alpha\in \Psi}m_\alpha\alpha,$$ with $\sum_{\alpha\in\Psi}n_\alpha=r$ and $\sum_{\alpha\in\Psi}m_\alpha= s+1$. Moreover since  $V\cong_{\lie g} V(\nu)$, we have
$$\dim \Hom_{\lie g}(V, P(\mu,r)[s]) = [P(\mu,r):V(\nu,s)]=\dim\Hom_{\lie g}(V(\nu), S^{s-r}(\lie g_{\ad})\otimes V(\mu))\ne 0,$$  which by~\lemref{PRE10}\eqref{PRE10.iv}
implies that
\begin{gather*}\nu=\mu-\sum_{\beta\in R} k_\beta\beta,\qquad k_\beta\in\bz_+, \quad \sum_{\beta\in R} k_\beta\le s-r,
\intertext{which gives}
\lambda-\nu=\sum_{\alpha\in\Psi}n_\alpha\alpha+\sum_{\beta\in R}k_\beta\beta.
\end{gather*}

If $\zeta=\nu$, we now get by using~\lemref{CAT500} that
$$\sum_{\alpha\in\Psi} m_\alpha= s+1\le \sum_{\alpha\in\Psi}n_\alpha+\sum_{\beta\in R} k_\beta\le s$$
which is absurd and hence we must have $\zeta-\nu=\beta_0$ for some $\beta_0\in R$. This gives
$$\lambda-\zeta= \sum_{\alpha\in \Psi}m_\alpha\alpha=\sum_{\alpha\in\Psi}n_\alpha\alpha+
\sum_{\beta\in R}k_\beta\beta+\beta_0 , $$and hence $$\sum_{\alpha\in\Psi}m_\alpha=
\sum_{\alpha\in\Psi}n_\alpha+\sum_{\beta\in R}k_\beta+1 \le s+1.$$ By Lemma 2.3
 again, this means that $k_\beta=0$ for all $\beta\notin\Psi$ and $\beta_0\in\Psi$.
  Therefore $\nu\le_\Psi\zeta\le_\Psi\lambda$ and  $d_\Psi(\nu,\lambda)=s$,
   hence $(\nu,s)\in\Gamma$ which is a contradiction.
\end{pf}

\subsection{}\label{CAT600}
In the following Proposition we gather the properties of projectives in~$\cal G_2[\le_\Psi\lambda]$ that will be needed later.
\begin{prop}
Let~$\lambda'\le_\Psi\lambda\in P^+$. Let  $\Gamma=\Lambda(\le_\Psi\lambda)$
or~$\Gamma=\Lambda([\lambda',\lambda]_\Psi)$.
\begin{enumerate}[{\rm(i)}]
 \item\label{CAT600.i} If $V\in\Ob\cal G_2[\Gamma]$
  is concentrated in degree~$r$,
  then $P(V)^\Gamma$ is the projective cover of~$V$ in~$\cal
  G_2[\Gamma]$. In particular, if $(\mu,r)\in\Gamma$, then
  $P(\mu,r)^\Gamma$ is the projective cover of~$V(\mu,r)$ in~$\cal G_2[\Gamma]$ and
   the category $\cal G_2[\Gamma]$ has enough projectives.
\item\label{CAT600.ia} Suppose that $(\xi,r)\notin\Gamma$ and
$$
\lambda-\xi=\sum_{\alpha\in R} n_\alpha\alpha,\qquad n_\alpha\in\bz_+,\,\sum_{\alpha\in R} n_\alpha\le r.
$$
Then $P(\xi,r)\in\Ob\cal G_2[\Lambda\setminus\Gamma]$.
\item\label{CAT600.ii} For
$(\mu,r),(\nu,s)\in\Gamma$, we have
$$[P(\mu,r)^\Gamma:V(\nu,s)]=\dim\Hom_{\cal
G_2}(P(\nu,s)^\Gamma,P(\mu,r)^\Gamma)\ne 0,$$  only if
~$\nu\le_\Psi\mu$ and $d_\Psi(\nu,\mu)=s-r$.
 \item\label{CAT600.iv} $\Hom_{\cal
G_2}(P(\nu,s),P(\mu,r))\cong\Hom_{\cal
G_2}(P(\nu,s)^\Gamma,P(\mu,r)^\Gamma)$ and this isomorphism is
compatible with compositions.
\end{enumerate}
\end{prop}
\begin{pf}
It follows from~\lemref{CAT530} that if $M\in\Ob\cal G_2[\Gamma]$
and $(\mu,r)\in\Gamma$, then
$$
\Hom_{\cal G_2}(P(\mu,r)_{\Lambda\setminus\Gamma},M)=0,
$$
or equivalently that
\begin{equation}\label{CAT600.10}
\Hom_{\cal G_2}(P(\mu,r)^\Gamma,M)\cong\Hom_{\cal
G_2}(P(\mu,r),M).
\end{equation}
In particular, $\Hom_{\cal G_2}(P(\mu,r)^\Gamma,-)$ is exact on~$\cal G_2[\Gamma]$
and hence~$P(\mu,r)^\Gamma$ is projective. The fact that it is isomorphic to the projective cover of~$V(\mu,r)$
in~$\cal G_2[\Gamma]$ is immediate.

To prove~\eqref{CAT600.ia}, suppose that~$[P(\xi,r):V(\nu,s)]\not=0$ for some~$(\nu,s)\in\Gamma$.
By~\eqref{CAT40.30}
$$
[P(\xi,r):V(\nu,s)]=\dim\Hom_{\lie g}(V(\nu),S^{s-r}(\lie g_{\ad})\tensor V(\xi)),
$$
hence by~\lemref{PRE10}\eqref{PRE10.iv}
$$
\xi-\nu=\sum_{\alpha\in R} m_\alpha\alpha,\qquad m_\alpha\in\bz_+,\,\sum_{\alpha\in R}m_\alpha\le s-r.
$$
Using~\lemref{CAT500} we conclude that $\nu\le_\Psi\xi\le_\Psi\lambda$,
and $d_\Psi(\xi,\lambda)=r$. This forces~$(\xi,r)\in\Gamma$ which is
a contradiction.

To prove part~\eqref{CAT600.ii}, note that~\eqref{CAT40.30} and~\lemref{CAT530} imply that
$$
[P(\mu,r)^\Gamma:V(\nu,s)]=[P(\mu,r): V(\nu,s)]=\dim\Hom_{\lie
g}(V(\nu),S^{s-r}(\lie g_{\ad})\tensor V(\mu)).
$$
In particular, if $[P(\mu,r)^\Gamma:V(\nu,s)]\ne 0$, we must have $$\mu-\nu=\sum_{\alpha\in R} k_\alpha\alpha,\qquad k_\alpha\in\bz_+,$$ 
and~$0\le\sum_{\alpha\in R} k_\alpha\le s-r$. 
Since $(\mu,r)\in\Gamma$ we write $\lambda-\mu=\sum_{\alpha\in\Psi}\ell_\alpha\alpha$ for some $\ell_\alpha\in\bz_+$ with $\sum_{\alpha\in\Psi}\ell_\alpha=r$. Then $$\lambda-\nu=\sum_{\alpha\in\Psi}\ell_\alpha\alpha+\sum_{\alpha\in R} k_\alpha\alpha,\qquad \sum_{\alpha\in \Psi} \ell_\alpha+
\sum_{\alpha\in R} k_\alpha\le s.$$ Since $(\nu,s)\in\Gamma$ it follows 
from \lemref{CAT500} that
$$
\sum_{\alpha\in\Psi} \ell_\alpha+\sum_{\alpha\in R} k_\alpha=s,
$$
hence~$k_\alpha=0$ if~$\alpha\notin\Psi$. This implies that~$\mu\le_\Psi\nu$ and
$$
\sum_{\alpha\in R} k_\alpha=\sum_{\alpha\in\Psi} k_\alpha=s-r,
$$
hence $d_\Psi(\mu,\nu)=s-r$.  The equality in~\eqref{CAT600.ii} is a standard property of projectives.

To prove part~\eqref{CAT600.iv}, note that~\eqref{CAT600.10} implies
$$
\Hom_{\cal G_2}(P(\nu,s)^\Gamma,P(\mu,r)^\Gamma)\cong\Hom_{\cal G_2}(P(\nu,s),P(\mu,r)^\Gamma).
$$
 Since
$\Hom_{\cal G_2}(P(\nu,s),P(\mu,r))$ maps onto
$\Hom_{\cal G_2}(P(\nu,s),P(\mu,r)^\Gamma)$, part~\eqref{CAT600.iv} now follows
by using~\propref{CAT40} and~\lemref{CAT530}, which give
\begin{equation*}
\dim\Hom_{\cal G_2}(P(\nu,s)^\Gamma,P(\mu,r)^\Gamma)=
\dim\Hom_{\cal G_2}(P(\nu,s),P(\mu,r)).\qedhere
\end{equation*}
\end{pf}

\subsection{}\label{CAT640} We can now prove one part of~\thmref{equivalence}.
 \begin{prop} Assume that $\Gamma=\Lambda(\le_\Psi\lambda)$.
The category~$\cal G_2[\Gamma]$ is equivalent to the category  of right finite dimensional modules for the
associative algebra $\End_{\cal G_2[\Gamma]}P(\Gamma)$, where $P(\Gamma)=\bigoplus_{(\mu,s)\in\Gamma} P(\mu,s)^{\Gamma}$.  Moreover, if we set
\begin{equation}\label{CAT640.10}
(\End_{\cal G_2[\Gamma]}P(\Gamma))[k]=\bigoplus_{(\mu,r),(\nu,s)\in\Gamma: r-s=k}\Hom_{\cal G_2[\Gamma]}(P(\mu,r)^{\Gamma}, P(\nu,s)^\Gamma),
\end{equation}
then we have an isomorphism of $\bz_+$-graded associative algebras $$\End_{\cal G_2[\Gamma]}P(\Gamma)\cong \bs_\Psi^{\lie g}(\le_\Psi\lambda)^{op}.$$
In particular, $\cal G_2[\Gamma]$ is equivalent to the category of left finite dimensional $\bs_\Psi^{\lie g}(\le_\Psi\lambda)$-modules.

The corresponding results also hold for $\lambda'\le_\Psi\lambda$, $\Gamma=\Lambda([\lambda',\lambda]_\Psi)$
and $\bs_\Psi^{\lie g}([\lambda',\lambda]_\Psi)$.
\end{prop}
\begin{pf}
By~\propref{CAT600}\eqref{CAT600.i}, $P(\Gamma)$ is a projective generator of~$\cal G_2[\Gamma]$.
The equivalence of categories is then standard and is provided by the extact functor $\Hom_{\cal G_2}(P(\Gamma),-)$. To prove the second assertion,
note that by~\propref{CAT600}\eqref{CAT600.ii}, \eqref{CAT640.10} defines a grading on~$\End_{\cal G_2}P(\Gamma)$. Then by~\propref{CAT600}\eqref{CAT600.iv}
we have
$$
\End_{\cal G_2} P(\Gamma)\cong \bigoplus_{(\mu,r),(\nu,s)\in \Gamma} \Hom_{\cal G_2}(P(\mu,r),P(\nu,s))
$$
as $\bz_+$-graded associative algebras. Furthermore, observe that there exists a canonical isomorphism
$$\Hom_{\cal G_2}(P(\mu,r), P(\nu,s))\cong \Hom_{\lie g}(V(\mu), S^{r-s}(\lie g_{\ad})\otimes V(\nu))$$ given by restriction and moreover this map is compatible with compositions, in the sense that if $f\in \Hom_{\cal G_2}(P(\mu,r), P(\nu,s))$, $g\in\Hom_{\cal G_2}(P(\nu,s), P(\xi,k))$, then the following diagram commutes
$$
\begin{diagram}
\node{}\node{S^{r-s}(\lie g)\tensor V(\nu)}\arrow{se,l}{1\tensor g|_{1\tensor V(\nu,s)}}\\
\node{V(\mu)}\arrow{ne,l}{f|_{1\tensor V(\mu,r)}}\arrow{se,l}{(g\circ f)|_{1\tensor V(\mu,r)}}\node{}
\node{S^{r-s}(\lie g)\tensor S^{s-k}(\lie g)\tensor V(\nu)}\arrow{sw,l}{m_{S(\lie g)}\tensor 1}\\
\node{}\node{S^{r-k}(\lie g)\tensor V(\xi)},
\end{diagram}
$$
where~$m_{S(\lie g)}:S(\lie g)\tensor S(\lie g)\to S(\lie g)$ is the multiplication map.
The proposition now follows from \lemref{PRE10}\eqref{PRE10.v}  along with the observation that
the multiplication of ~$\bs^{\lie g}$ when restricted to~$1_\mu\bs^{\lie g}[r-s]1_\nu \tensor
1_\nu\bs^{\lie g}[s-k]1_\xi$ is just the natural map
$$
(V(\mu)^*\tensor S^{r-s}(\lie g)\tensor V(\nu))^{\lie g}\tensor (V(\nu)^*\tensor S^{s-k}(\lie g)\tensor V(\xi))\to
(V(\mu)^*\tensor S^{r-s}(\lie g)\tensor S^{s-k}(\lie g)\tensor V(\xi))^{\lie g}
$$
composed with~$1\tensor m_{S(\lie g)}\tensor 1$.
\end{pf}

\section{Homological properties of category~\texorpdfstring{$\cal G_2$}{G\_2}}\label{HOM}

 \subsection{}
\label{CAT60} Given $(\mu,r)\in\Lambda$ and $0\le j\le\dim\lie
g$, set
$$
P_j(\mu,r)= S(\lie g_{\ad})\otimes (\bigwedge^j\lie g_{\ad}\otimes
V(\mu))[j+r],
$$ where we regard $(\bigwedge^j\lie g_{\ad}\otimes V(\mu))[j+r]$ as a
$\lie g\ltimes\lie g_{\ad}$-module concentrated in degree~$(j+r)$
with the canonical $\lie g$-action. In particular,
$P_0(\mu,r)=P(\mu,r)$ and the modules~$P_j(\mu,r)$ are projective
in~$\cal G_2$.

For $j\in\bz_+$ with $j>0$, let $d_j: P_j(\mu,r)\to
P_{j-1}(\mu,r)$ be the linear map obtained by extending
$$d_j(u\otimes (x_1\wedge x_2\wedge\cdots\wedge x_j)\tensor v)=
\sum_{s=1}^j (-1)^{s-1} ux_s\otimes (x_1\wedge\cdots\wedge
x_{s-1}\wedge x_{s+1}\wedge\cdots\wedge x_j)\tensor v,$$ and
let~$d_0: P_0(\mu,r)\to V(\mu,r)$ be
$$
d_0(u\tensor v)=u v
$$
where $u\in S(\lie g_{\ad})$, $v\in V(\mu,r)$ and $x_s\in \lie g$
for $1\le s\le j$.
\begin{prop} Let $(\mu,r)\in\Lambda$, $N=\dim{\lie g}$.
The  sequence
$$
0\longrightarrow P_{N}(\mu,r)\xrightarrow{d_N} P_{N-1}(\mu,r)
\stackrel{d_{N-1}}\longrightarrow\cdots
\stackrel{d_2}{\longrightarrow}
P_{1}(\mu,r)\stackrel{d_1}\longrightarrow
P(\mu,r)\stackrel{d_0}{\longrightarrow} V(\mu,r)\longrightarrow 0
$$
is a projective resolution of~$V(\mu,r)$ in the category~$\cal
G_2$.
\end{prop}
\begin{pf} To prove the exactness, note that the above sequence is obtained by tensoring the Koszul complex for~$S(\lie g)$ with~$V(\mu,r)$ and introducing an appropriate
grading. The map~$d_j$ is just $D\tensor 1$ where~$D$ is the
Koszul differential. The exactness is then standard (see for
example~\cite{GS}). It is straightforward to check that $d_j$ is a
morphism in $\cal G_2$ for all $j$, which proves the proposition.
\end{pf}
\begin{cor} Let $j\ge 0$.  We have
$$d_j(u\otimes a\tensor
v)= (u\otimes 1)d_j(1\otimes a \tensor v),\qquad \forall\, u\in S(\lie g_{\ad}),\, a\in \bigwedge^j \lie g_{\ad},\, v\in V(\mu,r).$$
In particular $\Im d_j[s]=0$ if $j+r>s$ and $$\Im
d_{s-r}[s]\cong_{\lie g} P_{s-r}(\mu,r)[s] \cong_{\lie g} \bigwedge^{s-r}\lie g_{\ad}\tensor V(\mu).$$
\end{cor}
\begin{pf}
The first two assertions are immediate. To prove the last assertion, observe that
since $\ker d_{s-r}[s]=0=\Im d_{s-r+1}[s]$,  we get $d_{s-r}(P_{s-r}(\mu,r)[s])\cong_{\lie g} P_{s-r}(\mu,r)[s]$.
\end{pf}

\subsection{}\label{CAT80} We can now compute the $\Ext$ spaces for all simple
objects in the category~$\cal G_2$.
\begin{prop}
For all~$j\ge 0$, $(\mu,r),(\nu,s)\in\Lambda$,
$$
\Ext^j_{{\cal G}_2}(V(\mu,r),V(\nu,s))\cong \begin{cases}
\Hom_{\lie g}(\textstyle\bigwedge^{j} \lie g_{\ad}\tensor
V(\mu),V(\nu)),&j=s-r, \\ 0,&
{\text{otherwise}.}
\end{cases}
$$
\end{prop}
\begin{pf}
For $j\ge 1$,  the short exact sequence
$$
0\to \Im d_j\to P_{j-1}(\mu,r)\to \Im d_{j-1}\to 0,
$$
yields the isomorphism,
$$
\Ext^j_{\cal G_2}(V(\mu,r),V(\nu,s))\cong \Ext^1_{\cal
G_2}(\Im d_{j-1},V(\nu,s)),
$$ and also the exact sequence
\begin{multline*}
0\to \Hom_{\cal G_2}(\Im d_{j-1},V(\nu,s))\to \Hom_{\cal G_2}(P_{j-1}(\mu,r),V(\nu,s))\to\\
\Hom_{\cal G_2}(\Im d_j,V(\nu,s))\to \Ext^1_{\cal G_2}(\Im
d_{j-1},V(\nu,s))\to 0.
\end{multline*}
We claim that $\Hom_{\cal G_2}(\Im d_j,V(\nu,s))=0$
unless~$j=s-r$ which proves that  $$\Ext^j_{\cal
G_2}(V(\mu,r),V(\nu,s))\cong \Ext^1_{\cal G_2}(\Im
d_{j-1},V(\nu,s))=0,\ \ j\ne s-r.$$

To prove the claim, note that by~\corref{CAT60} we may assume
$j+r\le s$ and
also that if  $v\in\Im
 d_j[s]$, then  $$v=\sum_p (u_p\otimes 1) d_j(1\otimes w_p),$$
where $u_p\in
S^{s-r-j}(\lie g_{\ad})$ and $w_p\in (\bigwedge^j \lie
g_{\ad}\otimes V(\mu))[j+r]$. Suppose first that $j+r<s$. Let $f\in \Hom_{\cal G_2}(\Im
d_j,V(\nu,s))$. Then
 $f(1\otimes w_p)=0$ and so~$f(v)=\sum_p u_p f(w_p)=0$, that is $f=0$.

Suppose that $j=s-r$. Since~$P_{s-r-1}(\mu,r)$ is the projective
cover of a semi-simple $\lie g\ltimes \lie g_{\ad}$-module concentrated in degree~$s-1$,
it follows that
$$\Hom_{\cal G_2}(P_{s-r-1}(\mu,r),V(\nu,s))=0,
$$
and hence
\begin{align*}
\Ext^{s-r}_{\cal G_2}(V(\mu,r),V(\nu,s))&\cong \Ext^1_{\cal
G_2}(\Im d_{s-r-1},V(\nu,s))\cong \Hom_{\cal G_2}(\Im
d_{s-r},V(\nu,s))\\&\cong \Hom_{\lie g}(\Im
d_{s-r}[s],V(\lambda)).
\end{align*}
The result follows by applying~\corref{CAT60}.
\end{pf}

\subsection{}\label{CAT900}

\begin{prop}
Let~$\lambda'\le_\Psi\lambda\in P^+$. Assume that
$\Gamma=\Lambda(\le_\Psi\lambda)$ or $\Gamma=\Lambda([\lambda',\lambda]_\Psi)$
Let~$(\mu,r)\in\Gamma$.
Then the induced map $d_j^\Gamma: P_j(\mu,r)^\Gamma\to  P_j(\mu,r-1)^{\Gamma}$ is a morphism of objects in $\cal G_2[\Gamma]$ and
the  sequence
$$
0\longrightarrow P_{N}(\mu,r)^\Gamma\xrightarrow{d_N^\Gamma} P_{N-1}(\mu,r)^\Gamma
\stackrel{d_{N-1}^\Gamma}\longrightarrow\cdots
\stackrel{d_2^\Gamma}{\longrightarrow}
P_{1}(\mu,r)^\Gamma\stackrel{d_1^\Gamma}\longrightarrow
P(\mu,r)^\Gamma\stackrel{d_0^\Gamma}{\longrightarrow} V(\mu,r)\longrightarrow 0
$$
is a projective resolution of~$V(\mu,r)$ in the category~$\cal
G_2[\Gamma]$. In particular for all~$(\nu,s)\in\Gamma$  we have
\begin{gather*}
\Ext^j_{\cal G_2[\Gamma]}(V(\mu,r),V(\nu,s))\cong\Ext^j_{\cal G_2}(V(\mu,r),V(\nu,s))\end{gather*}
and hence $$\Ext^j_{\cal G_2[\Gamma]}(V(\mu,r),V(\nu,s))= \begin{cases}
\Hom_{\lie g}(\bigwedge^j \lie g_{\ad}\tensor V(\mu),V(\nu)),& \nu\le_\Psi\mu,\,j=d_\Psi(\nu,\mu)=s-r\\
0,&\text{otherwise.}
                                                                                 \end{cases}$$
\end{prop}

\begin{pf}
Set~$$W_j(\mu,r)=(\bigwedge^j \lie g_{\ad}\tensor V(\mu))[r+j].$$
Since~$W_j(\mu,r)\in\Ob\cal G_2$ is concentrated in degree $r+j$
and is $\lie g$-semisimple, it is also a  semi-simple object of
$\cal G_2$ and so, we have a $\lie g\ltimes\lie g_{\ad}$
decomposition,
$$W_j(\mu,r)=W_j(\mu,r)_{\Lambda\setminus\Gamma}\oplus
W_j(\mu,r)^\Gamma$$  and hence  by
 Proposition~\ref{CAT40}\eqref{CAT40.iv} we have$$P_j(\mu,r)\cong_{\cal G_2} P(W_j(\mu,r))\cong_{\cal G_2}
P(W_j(\mu,r)_{\Lambda\setminus\Gamma})\oplus
P(W_j(\mu,r)^\Gamma).$$
Note that  $[W_j(\mu,r):V(\xi,j+r)]\not=0$ implies that
$$
\mu-\xi =\sum_{\alpha\in R} m_\alpha\alpha,\qquad
m_\alpha\in\bz_+, \quad\sum_{\alpha\in R} m_\alpha\le j,
$$
and so
\begin{equation*}
\lambda-\xi=\sum_{\alpha\in R} m'_\alpha\alpha,\qquad m_\alpha'\in\bz_+,\quad
\sum_{\alpha\in R} m'_\alpha\le r+j.\end{equation*}
It now follows from~\propref{CAT600}\eqref{CAT600.ia} that
$P(W_j(\mu,r)_{\Lambda\setminus\Gamma})\in\Ob\cal G_2[\Lambda\setminus\Gamma]$, hence
$$
P_j(\mu,r)_{\Lambda\setminus\Gamma}\cong_{\cal G_2} P(W_j(\mu,r)_{\Lambda\setminus\Gamma})
\oplus P(W_j(\mu,r)^\Gamma)_{\Lambda\setminus\Gamma}.
$$
By~\lemref{CAT530}, $P_j(\mu,r)_{\Lambda\setminus\Gamma}\in\Ob\cal G_2[\Lambda\setminus\Gamma]$
and so~$P_j(\mu,r)^\Gamma$ is a projective object in~$\cal G_2[\Gamma]$ and is isomorphic to
$P(W_j(\mu,r)^\Gamma)^\Gamma$. All assertions are now straightforward.
\end{pf}

\subsection{}\label{CAT95}
Let~$\lambda_\Psi=\sum_{\beta\in\Psi}\beta$.

\begin{lem}
\begin{enumerate}[{\rm(i)}]
\item\label{CAT95.i} If~$\mu,\mu+\lambda_\Psi\in P^+$, $\dim\Hom_{\lie g}(V(\mu+\lambda_\Psi),\bigwedge^{|\Psi|}\lie g_{\ad}\tensor V(\mu))\le 1$.
\item\label{CAT95.ii} There exists~$\mu\in P^+$ such that~$\mu+\lambda_\Psi\in P^+$ and
$\dim\Hom_{\lie g}(V(\mu+\lambda_\Psi),\bigwedge^{|\Psi|}\lie g_{\ad}\tensor V(\mu))=1$.
\end{enumerate}
\end{lem}
\begin{pf}
Using ~\lemref{CAT500} we see that  $(\bigwedge^{|\Psi|} \lie g_{\ad})_{\lambda_\Psi}=\bc\bigwedge_{\beta\in\Psi} x^+_\beta$ and hence $\dim (\bigwedge^{|\Psi|} \lie g_{\ad})_{\lambda_\Psi}=1$. Part~\eqref{CAT95.i} is now immediate
from~\lemref{PRE10}\eqref{PRE10.iv}. To prove~\eqref{CAT95.ii}, choose $k_i\in\bz_+$, $i\in I$ such that $$(\ad x_{\alpha_i}^+)^{k_i+1} \bigwedge_{\beta\in\Psi} x^+_\beta=0 = (\ad x_{\alpha_i}^-)^{k_i+1+\lambda_\Psi(h_i)} \bigwedge_{\beta\in\Psi} x^+_\beta,\quad\lambda_\Psi+\sum_{i\in I} k_i\omega_i\in P^+$$
and set~$\mu=\sum_{i\in I} k_i\omega_i$.
 It follows from \lemref{PRE10}\eqref{PRE10.iv} that
\begin{equation*}
\Hom_{\lie g}(V(\mu+\lambda_\Psi),\bigwedge^{|\Psi|}\lie g_{\ad}\tensor V(\mu))\cong \bc \bigwedge_{\beta\in\Psi} x_\beta^+.\qedhere
\end{equation*}
\end{pf}

\subsection{}\label{CAT100} Recall that for an abelian category $\cal C$ which has enough injectives or projectives, the global dimension, written $\gldim \cal C$, equals the minimal~$j$ such that~$\Ext^j_{\cal C}(M,N)=0$ for all~$M,N\in\Ob\cal C$.
Note that categories~$\cal G_2[\le_\Psi\lambda]$, $\lambda\in P^+$ and $\cal G_2[[\lambda',\lambda]_\Psi]$, $\lambda'\le_\Psi\lambda$ have enough projectives by~\propref{CAT530}.
We can now prove
\begin{thm}\label{thmCAT100} Let $\lambda'\le_\Psi\lambda\in P^+$. We have
 \begin{gather*}
\gldim\cal G_2[[\lambda',\lambda]_\Psi], \gldim \cal G_2[\le_\Psi\lambda]\le |\Psi|
\end{gather*}
and the upper bound is attained for some~$\lambda'\le_\Psi\lambda\in P^+$.
In particular, the algebras $\bs_{\Psi}^{\lie g}(\le_\Psi\lambda)$,
$\bs_{\Psi}^{\lie g}([\lambda',\lambda]_\Psi)$ have global dimension at most~$|\Psi|$ and
the upper bound is attained for some~$\lambda'\le_\Psi\lambda\in P^+$.
\end{thm}
\begin{pf}
Let~$\Gamma=\Lambda(\le_\Psi\lambda)$.
Since all objects in~$\cal G_2[\Gamma]$ have finite length,
to establish an upper bound  on the global dimension of $\cal G_2[\Gamma]$, it suffices to prove that
for all $(\mu,r), (\nu,k)\in\Gamma$ we have
$$\Ext^j_{\cal G_2[\Gamma]}(V(\mu,r),V(\nu,k))=0,\qquad\forall\,j>|\Psi|.
$$
By~\propref{CAT80}, this amounts to proving that for $\nu\le_\Psi\mu$, $j=d_\Psi(\nu,\mu)$
$$
\Hom_{\lie g}(\bigwedge^{j} \lie g_{\ad}\tensor V(\mu),V(\nu))
\ne 0\implies |\Psi|\ge j.
$$

A weight  of $\bigwedge^j\lie g$ is the sum of at most $j$ distinct roots.
Hence by \lemref{PRE10}\eqref{PRE10.iv} it follows that
\begin{equation}\label{AP40.20}
\Hom_{\lie g}(\bigwedge^{j} \lie g_{\ad}\tensor V(\mu),V(\nu))\not=0
\end{equation}
only if $$\mu-\nu =\sum_{\alpha\in R} m_\alpha\alpha,\qquad m_\alpha\in\{0,1\}, \quad \sum_{\alpha\in R}m_\alpha \le j.
$$
Applying~\lemref{CAT500}, we immediately conclude that~$m_\alpha=0$ unless~$\alpha\in\Psi$ and also that $\sum_{\alpha\in \Psi}m_\alpha=j$.
Since~$m_\alpha\in\{0, 1\}$, it follows that $|\Psi|\ge j$.

To prove that the upper bound is attained, observe that
by~\lemref{CAT95}\eqref{CAT95.ii}, there exists~$\mu\in P^+$ such that~$\mu+\lambda_\Psi\in P^+$ and $\Hom_{\lie g}(V(\mu+\lambda_\Psi),\bigwedge^{|\Psi|}\lie g_{\ad}\tensor V(\mu))\not=0$. Then
$(\mu,|\Psi|)\in\Lambda(\mu+\lambda_\Psi,\Psi)$ and
we have
\begin{align*}
\dim\Ext^{|\Psi|}_{\cal G_2}(V(\mu+\lambda_\Psi,0),V(\mu,|\Psi|))&=\dim\Hom_{\lie g}(\bigwedge^{|\Psi|} \lie g_{\ad}
\tensor V(\mu+\lambda_\Psi),V(\mu))\\
&=\dim\Hom_{\lie g}(V(\mu+\lambda_\Psi),\bigwedge^{|\Psi|}\lie g_{\ad}\tensor V(\mu))=1
.\qedhere
\end{align*}
\end{pf}
\subsection{} \begin{prop} The algebras~$\bs_\Psi^{\lie g}$, ~$\bs_\Psi^{\lie g}(\lambda\le_\Psi)$ have left global dimension~$|\Psi|$.
\end{prop}
\begin{pf}
Given~$\mu\in P^+$,
let~$S_\mu$ be the left simple $\bs_\Psi^{\lie g}$-module $S_\mu$ corresponding to the idempotent~$1_\mu$.
Then its projective cover in the category $\bs^{\lie g}_\Psi-\mof$ of finite dimensional left $\bs^{\lie g}_\Psi$-modules is
$$
P_\mu=\bs^{\lie g}_\Psi 1_\mu=
\bigoplus_{\nu\le_\Psi\mu} 1_\nu\bs_\Psi^{\lie g} 1_\mu=\bs_\Psi^{\lie g}(\le_\Psi\mu) 1_\mu.
$$
Similarly, if~$\nu\le_\Psi\mu$, $P_\nu=\bs_\Psi^{\lie g}(\le_\Psi\mu)1_\nu$.
In particular, $[P_\mu:S_\nu]=0$ unless~$\nu\le_\Psi\mu$. Proceeding by induction we conclude that~$S_\mu$
has a projective resolution in $\bs_\Psi^{\lie g}-\mof$
$$
\cdots \to  P^1_\mu\to  P_\mu\to S_\mu\to 0
$$
in which~$[ P^j_\mu:  S_\xi]=0$ unless~$\xi\le_\Psi\mu$.
Therefore, this projective resolution can be regarded as
a projective resolution for~$S_\mu$ in the category $\bs_\Psi^{\lie g}(\le_\Psi\mu)-\mof$, which is
equivalent to the category~$\cal G_2[\le_\Psi\mu]$ by~\propref{CAT640}, and~$\Ext^j_{\bs^{\lie g}_\Psi-\mof}
(S_\mu,S_\nu)=0$ unless~$\nu\le_\Psi\mu$. Thus, we have
\begin{multline*}
\dim\Ext^j_{\bs_\Psi^{\lie g}-\mof}(S_\mu,
S_\nu)=\dim\Ext^j_{\bs_\Psi^{\lie g}(\le_\Psi\mu)-\mof} (S_\mu,S_\nu)\\=\dim\Ext^j_{\cal G_2[\le_\Psi\mu]}(V(\mu,0),V(\nu,d_\Psi(\nu,\mu))),
\end{multline*}
and the result follows from~\thmref{thmCAT100}.

The argument for~$\bs^{\lie g}_\Psi(\lambda\le_\Psi)$ is similar, with~$\bs^{\lie g}_\Psi(\le_\Psi\mu)$ replaced by~$\bs^{\lie g}_\Psi([\lambda,\mu]_\Psi)$.\end{pf}
\subsection{}\label{CAT105}
Let~$A=\bigoplus_{r\in\bz_+} A[r]$ be a $\bz_+$-graded associative algebra such that~$A[0]$ is semi-simple. We can regard
$A[0]$ as a graded left $A$-module concentrated in degree~$0$ via the canonical projection $A\to A/\bigoplus_{j>0} A[j]\cong A[0]$.
Following~\cite[Definition~1.2.1]{BGS}, a given grading on~$A$ is said to be Koszul
if~$A[0]$ admits a projective resolution
$$
\cdots\to P^2\to P^1\to P^0\to A[0]\to 0
$$
in the category of left $\bz_+$-graded $A$-modules such that $P^r$ is generated by~$P^r[r]$ as a graded $A$-module.

Suppose that~$\dim A[r]<\infty$ for all~$r\in\bz_+$ and~$A[0]$ is commutative, so that in particular we have~$A[0]=\bigoplus_{\mu\in F} \bc 1_\mu$,
where the~$1_\mu$ are pairwise orthogonal idempotents and~$F$ is a finite index set. Then
we have the following numerical criterion for Koszulity (\cite[Theorem~2.11.1]{BGS}). Let~$S_\mu=\bc 1_\mu$ be the simple left $A$-module corresponding to~$1_\mu$.
Set
$$
H( A ,t)=(H( A ,t)_{\mu,{\nu}})_{\mu,{\nu}\in F}, \quad
H(E( A ),t)=(H(E( A ),t)_{\mu,{\nu}})_{\mu,{\nu}\in F}
$$
where
\begin{align*}
&H( A ,t)_{\mu,{\nu}}=\sum_{i\ge 0} t^i \dim(1_\mu  A  [i] 1_{\nu}), \\
&H(E( A  ),t)_{\mu,{\nu}}=\sum_{i\ge 0} t^i \dim \Ext^i_{ A -\mof} (S_\nu,S_{\mu})
\end{align*}
are formal power series in~$t$ with coefficients in~$\bz_+$. The matrix~$H(A,t)$ is the Hilbert matrix of~$A$, while $H(E(A),t)$ is the
Hilbert matrix of the Yoneda algebra of~$A$. The algebra~$A$ is Koszul if and only if $H(E(A),-t)$ is the inverse of~$H(A,t)$.

\subsection{}\label{CAT110} The following proposition completes the proof of~\thmref{equivalence}.

\begin{prop}
Let~$\lambda'\le_\Psi\lambda\in P^+$.
The natural grading on $\bs_\Psi^{\lie g}(\le_\Psi\lambda)$ and $\bs_\Psi^{\lie g}([\lambda',\lambda]_\Psi)$ is
Koszul.
\end{prop}
\begin{pf}
Let~$F=(\le_\Psi\lambda)$ or~$F=[\lambda',\lambda]_\Psi$ and set $A =\bs_\Psi^{\lie g}(F)$.
For~$\mu\in F$ let $S_\mu$ be the simple left $ A $-module corresponding to the
idempotent~$1_\mu$ and let~$P_\mu= A  1_\mu$ be its projective cover. Then~$[P_\mu:S_\nu]=\dim(1_\nu A 1_\mu)$.

By~\propref{CAT640}, the category $ A -\mof$ is equivalent to the category~$\cal G_2[F]$.
\propref{CAT600} gives
$$
\dim(1_\mu  A  [j] 1_{\nu})\not=0\,\implies\, {\mu}\le_\Psi\nu,\, j=d_\Psi({\mu},\nu).
$$
and using ~\propref{CAT900} we also have
$$
\Ext_{ A -\mof}^j (S_{\nu},S_{{\mu}})\not=0\,\implies\, {\mu}\le_\Psi\nu,\,j=d_\Psi({\mu},\nu).
$$
This shows that  the  matrices $H( A  ,t)$ and~$H(E(  A  ),t)$, and so their product,
are all  upper triangular.
Moreover, for ~${\mu}\le_\Psi\nu \in F$ we have
\begin{align*}
\sum_{\xi\in F} &H(E( A ) ,-t)_{\mu,\xi} H( A ,t)_{\xi,{\nu}}\\
&=\sum_{{\mu}\le_\Psi\xi\le_\Psi\nu} (-1)^{d_\Psi(\mu,{\xi})}
t^{d_\Psi(\xi,{\nu})+d_\Psi(\mu,\xi)}
[P_\nu:S_\xi]\dim\Ext^{d_\Psi(\mu,{\xi})}_{ A -\mof}(S_\xi,S_{\mu})\\
&=t^{d_\Psi({\mu},\nu)} \sum_{j\ge 0}\sum_{{\mu}\le_\Psi\xi\le_\Psi\nu} (-1)^j [P_\nu:S_\xi] \dim\Ext^j_{ A -\mof}(S_\xi,S_{\mu})\\
&=t^{d_\Psi({\mu},\nu)} \sum_{j\ge 0}(-1)^j\dim\Ext^j_{ A -\mof}(P_\nu,S_{\mu})\\
&=\delta_{\mu,{\nu}} t^{d_\Psi({\mu},\nu)}=\delta_{\mu,{\nu}}.
\end{align*}
Thus, the matrix $H(E(A),-t)$ is the inverse of the matrix~$H(A,t)$ and so ~$ A $ is Koszul by ~\cite[Theorem~2.11.1]{BGS}.
\end{pf}

\section{Koszulity of~\texorpdfstring{$\bs_\Psi^{\lie g}$}{\bfseries S\_Psi g}}

\subsection{}\label{SALG10} We shall use the following elementary result repeatedly.
\begin{lem} Let $M\in\Ob\cal F(\lie g)$.
 There exists an isomorphism $$ M\cong
\bigoplus_{\nu\in P^+}V(\nu)\tensor (V(\nu)^*\tensor M)^{\lie
g},$$ of $\lie g$-modules. In particular, if $M'\in\Ob\cal
F(\lie g)$, then we have $$(M'\otimes M)^{\lie
g}\cong\bigoplus_{\nu\in P^+} (M'\otimes V(\nu))^{\lie g}\otimes
(V(\nu)^*\tensor M)^{\lie g}.$$
\end{lem}
\begin{pf} It suffices to prove the Lemma in the case when $M$ is a direct sum of copies of $V(\xi)$ for some $\xi\in P^+$. In this case, it is clear that the map $V(\xi)\tensor \Hom_{\lie g}(V(\xi), M)\to M$ given by extending $v\tensor f\mapsto f(v)$, induces an isomorphism of $\lie g$-modules. The Lemma follows by noticing that $\Hom_{\lie g}(V(\xi), M)\cong (V(\xi)^*\tensor M)^{\lie g}$.
\end{pf}

\subsection{}\label{SALG20} Consider the  canonical surjection $\Pi_S: T(\lie g)\to S(\lie
g)$ and let $\boldsymbol{\Pi}_S: \bt\to\bs$ be the map $1\otimes\Pi_S\otimes 1$.
Identify the $\lie g$-submodule of~$T^2(\lie g)$ spanned by the elements
$\{x\otimes y-y\otimes x: x,y\in\lie g\}$ (respectively, by the elements
$\{x\otimes y+y\otimes x:x,y\in\lie g\}$) with $\bigwedge^2 \lie
g $ (respectively, with $S^2(\lie g)$). Then $\ker\Pi_S$ is the ideal generated by $\bigwedge^2\lie g$. Recall also that we have fixed $\psi\in P$ such that $\Psi=\Psi(\psi)$ is a subset of $R^+$.
We say that a subset $F$ of $P^+$ is interval closed if for $\mu,\nu\in F$ with $\mu\le_\Psi\nu$ we have $[\mu,\nu]_\Psi\subset F$.
It is clear that $\boldsymbol{\Pi}_S: \bt_\Psi(F)\to\bs_\Psi(F)$.
We first prove
\begin{lem}
\begin{enumerate}[{\rm(i)}]
\item\label{SALG20.i} The kernel of the restriction $\boldsymbol{\Pi}^{\lie g}_S:\bt^{\lie g}\to
\bs^{\lie g}$  is generated by $(\bv^{\circledast}\tensor \bigwedge^2 \lie
g\tensor \bv)^{\lie g}\subset \bt^{\lie g}[2]$.
\item\label{SALG20.ii} Let $F\subset P^+$ be interval closed.   The kernel of $\boldsymbol{\Pi}^{\lie g}_S$ restricted to $\bt_\Psi^{\lie g}(F)$  is generated by $(\bv^{\circledast}\tensor \bigwedge^2 \lie
g\tensor \bv)^{\lie g}\cap \bt_\Psi^{\lie g}(F)[2].$
\end{enumerate}

 Analogous
statements hold for $\be^{\lie g}$ and $\be^{\lie g}_\Psi(F)$.\end{lem}
\begin{pf}
 Note that for $k=0,1$ we have $$\bt^{\lie g}[k]=\bs^{\lie g}[k],$$ and
 hence $\ker\boldsymbol{\Pi}^{\lie g}_S[k]=0$ in these cases.
 For~$k\ge 2$, we have
$$
\ker\boldsymbol{\Pi}^{\lie g}_S\cap \bt^{\lie g}[k]=(\bv^{\circledast}\tensor
\ker\Pi_S[k]\tensor \bv)^{\lie g}=\sum_{j=0}^{k-2} ( \bv^{\circledast}\tensor
(T^j(\lie g)\tensor \bigwedge^2\lie g\tensor T^{k-j-2}(\lie
g))\tensor\bv)^{\lie g}.
$$
Part~\eqref{SALG20.i} of the Lemma follows if we  prove  for all ~$\lambda,\mu\in P^+$ and
for all  $j,k$ with  $k\ge 2$ and~$0\le j\le k-2$, that
$$
1_\lambda (\bv^{\circledast}\tensor (T^j(\lie g)\tensor \bigwedge^2 \lie
g\tensor T^{k-j-2}(\lie g))\tensor \bv)^{\lie g} 1_\mu \subset
1_\lambda \bt^{\lie g}[j] (\bv^{\circledast}\tensor \bigwedge^2 \lie
g\tensor\bv)^{\lie g} \bt^{\lie g}[k-j-2] 1_\mu,
$$
or  that
\begin{multline}\label{kernel}
(V(\lambda)^*\tensor (T^j(\lie g)\tensor \bigwedge^2 \lie g\tensor T^{k-j-2}(\lie g))\tensor V(\mu))^{\lie g}\\
=\sum_{\nu,\xi\in P^+} (V(\lambda)^* \tensor T^j(\lie g)\tensor
V(\nu))^{\lie g}(V(\nu)^*\tensor \bigwedge^2\lie g\tensor
V(\xi))^{\lie g} (V(\xi)^*\tensor T^{k-j-2}(\lie g)\tensor
V(\mu))^{\lie g}
\end{multline}
(the product in the right hand side is taken in~$\bt^{\lie g}$).
But this follows from~\lemref{SALG10} and the proof of part (i) is complete.

To prove~\eqref{SALG20.ii} suppose now that $\lambda\le_\Psi\mu$ and that $k=d_\Psi(\lambda,\mu)$. Again if $k=0,1$ there is nothing to prove.
If $k\ge 2$, then by~\lemref{PRE10} we see that
\begin{align*} (V(\lambda)^* \tensor T^j(\lie g)\tensor
V(\nu))^{\lie g}\ne 0&\implies\nu-\lambda=\sum_{\alpha\in R} m_\alpha\alpha,\qquad  m_\alpha\in\bz_+,\, \sum_{\alpha\in R} m_\alpha\le j,\\
(V(\nu)^*\tensor \bigwedge^2\lie g\tensor
V(\xi))^{\lie g}\ne 0&\implies\xi-\nu=\sum_{\alpha\in R}n_\alpha\alpha,\qquad n_\alpha\in\bz_+, \, \sum_{\alpha\in R}n_\alpha\le 2,\\
(V(\xi)^*\tensor T^{k-j-2}(\lie g)\tensor
V(\mu))^{\lie g}\ne 0 &\implies \mu-\xi=\sum_{\alpha\in R}\ell_\alpha\alpha,\qquad\ell_\alpha\in\bz_+,\, \sum_{\alpha\in R}\ell_\alpha\le k-j-2.
\end{align*}
But now~\lemref{CAT500} gives  $$\alpha\notin \Psi\implies m_\alpha= n_\alpha=\ell_\alpha=0,\qquad \sum_{\alpha\in\Psi}m_\alpha+n_\alpha+\ell_\alpha=j.
$$ 
This means that
$$
\lambda\le_\Psi\nu\le_\Psi\xi\le_\Psi\mu, \qquad d_\Psi(\lambda,\nu)=j,\qquad d_\Psi(\nu,\xi)=2,\qquad d_\Psi(\xi,\mu)=k-j-2,
$$
and since~$F$ is interval closed we get~$\nu,\xi\in F$.
In other words, we have proved that in the case when $\lambda\le_\Psi\mu$ and $d_\Psi(\lambda,\mu)=k$, the right hand side of \eqref{kernel} is in the ideal of $\bt^{\lie g}_\Psi(F)$ generated by $(\bv^{\circledast}\tensor \bigwedge^2 \lie
g\tensor \bv)\cap \bt_\Psi^{\lie g}(F)[2]$, which establishes part~\eqref{SALG20.ii} of the Lemma.
\end{pf}

\subsection{}\label{SALG0}
Before we continue proving the main result of this section, we  summarize for the reader's convenience some standard facts about $\bz_+$-graded associative algebras.  The details
can be found in \cite{BGS}.

Suppose that ~$A=\bigoplus_{r\in\bz_+} A[r]$ be a $\bz_+$-graded associative algebra and that
$A[0]$ is semi-simple. Clearly  $A[r]$ is an $A[0]$-bimodule for all $r\ge 0$ and we let ~$T^r_{A[0]}(A[1])$ be the $r$-fold tensor product of the $A[0]$-bimodule $A[1]$ over~$A[0]$.
Setting,$$
T^0_{A[0]}(A[1])=A[0],\qquad T_{A[0]}(A[1])=\bigoplus_{r\in\bz_+} T^r_{A[0]}(A[1]),
$$
we see that $T_{A[0]}(A[1])$ is a $\bz_+$-graded associative algebra and the assignment $\bom(a)=a$, $a\in A[r]$, $r=0,1$ extends to  a  canonical homomorphism $
\bom: T_{A[0]}(A[1])\to A
$
of $\bz_+$-graded associative algebras and $A[0]$-bimodules.
 The algebra $A$ is said to be {\em quadratic} if~$\bom$
is surjective and $\ker\bom$ is generated by $\ker\bom\cap T^2_{A[0]}(A[1])$.

 The Koszul complex for a quadratic algebra~$A$
is constructed as follows.
Set
\begin{gather*}
N^0_A=A[0],\qquad N^1_A=A[1],\\
N_A^r = \bigcap_{j=0}^{r-2} T^j_{A[0]}(A[1])\tensor_{A[0]} (\ker\bom\cap T^2_{A[0]}(A[1]))\tensor_{A[0]} T^{r-j-2}_{A[0]}(A[1]),\qquad r\ge 2.
\end{gather*}
Regard~$N^r$ as a graded left $A$-module concentrated in degree~$r$ via the canonical
projection $A\twoheadrightarrow A/\bigoplus_{j>0} A[j]\cong A[0]$.
Define a map
$$
\partial_r: A\tensor_{A[0]} N^r_A\to A\tensor_{A[0]} N^{r-1}_A
$$
by restricting the map
\begin{align*}
A\tensor_{A[0]} T^r_{A[0]}(A[1])&\to A\tensor_{A[0]} T^{r-1}_{A[0]}(A[1])\\
x\tensor a_1\tensor\cdots\tensor a_r&\mapsto xa_1\tensor a_2\cdots\tensor a_r,\qquad x\in A,\, a_j\in A[1],\, 1\le j\le r.
\end{align*}
Then $\partial_r$ is a homomorphism of graded left $A$-modules and of $A[0]$-bimodules and $\partial_{r-1}\partial_r=0$. The complex of $\bz_+$-graded left $A$-modules
$$
\cdots \stackrel{\partial_{r+1}}\longrightarrow A\tensor_{A[0]} N^r_A\stackrel{\partial_r}\longrightarrow A\tensor_{A[0]} N^{r-1}_A
\stackrel{\partial_{r-1}}{\longrightarrow}\cdots
\stackrel{\partial_2}\longrightarrow A\tensor_{A[0]} A[1]\stackrel{\partial_1}\longrightarrow A
$$
called the Koszul complex of~$A$. By~\cite[Theorem~2.6.1]{BGS}, a quadratic algebra is Koszul if and only if its Koszul complex is exact.

\subsection{}\label{SALG30} Let $F$ be a subset of $P^+$ and assume that $F$ is interval closed.
Recall that $$\bt^{\lie g}[0]=\bigoplus_{\lambda\in P^+} \bc 1_\lambda,\ \  \bt_\Psi^{\lie g}(F)[0]=\bigoplus_{\lambda\in F}\bc 1_\lambda,$$ is a semisimple, commutative algebra.
  Clearly, if ~$M,N$ are $\bt^{\lie g}[0]$-bimodules, then for any $\lambda\in P^+$ we have,$$M 1_\lambda\tensor_{\bt^{\lie g}[0]} 1_\lambda N=
M 1_\lambda\tensor 1_\lambda N.$$ Further, given a $\bt^{\lie g}_\Psi(F)[0]$-bimodule~$M$, we can always
regard~$M$ as a $\bt^{\lie g}[0]$-bimodule, by letting $1_\mu$, $\mu\notin F$ act trivially and in that case  we have
$$
M\tensor_{\bt^{\lie g}[0]} N=M\tensor_{\bt^{\lie g}_\Psi(F)[0]} N.
$$
In particular, this implies that $T_{\bt_\Psi^{\lie g}(F)[0]}(\bt_\Psi^{\lie g}(F)[1])$
is canonically isomorphic to a subalgebra of $T_{\bt^{\lie g}[0]}(\bt^{\lie g}[1])$
and the restriction of~$\bom$ maps it to $\bt_\Psi^{\lie g}(F)$.

Furthermore, we get
\begin{gather*}\bt^{\lie g}[1]\otimes_{\bt^{\lie g}[0]}\bt^{\lie g}[1]\cong\bigoplus_{\lambda,\mu,\nu\in P^+} (V(\lambda)^*\otimes\lie g\otimes V(\mu))^{\lie g}\otimes(V(\mu)^*\otimes\lie g\otimes V(\nu))^{\lie g},\\\bt_\Psi^{\lie g}(F)[1]\otimes_{\bt^{\lie g}[0]}\bt_{\Psi}^{\lie g}(F)[1]\cong\bigoplus_{(\lambda,\mu,\nu)\in\mathbf F} (V(\lambda)^*\otimes\lie g\otimes V(\mu))^{\lie g}\otimes(V(\mu)^*\otimes\lie g\otimes V(\nu))^{\lie g},\end{gather*} where $$\mathbf F=\{(\lambda,\mu,\nu)\in F^3:\lambda\le_\Psi\mu\le_\Psi\nu, \ d_{\Psi}(\lambda,\mu)=d_\Psi(\mu,\nu)=1\}.$$

 Applying~\lemref{SALG10} now gives,
\begin{align*}
\bt^{\lie g}[1]\otimes_{\bt^{\lie g}[0]}\bt^{\lie g}[1]&\cong\bigoplus_{\lambda,\nu\in P^+}
(V(\lambda)^*\otimes\lie g\otimes\lie g\otimes V(\nu))^{\lie g},\\ &\cong \bt^{\lie g}[2],
\end{align*}
while an argument identical to the one given in the proof of \lemref{SALG20}\eqref{SALG20.ii} shows that
$$\bt_\Psi^{\lie g}(F)[1]\otimes_{\bt_\Psi^{\lie g}[0]}\bt_\Psi^{\lie g}(F)[1]\cong \bt_\Psi^{\lie g}(F)[2].$$
More generally, we see now that we have an isomorphism of $\bt^{\lie g}[0]$-bimodules,
$$
(\bt^{\lie g}[1])_{\bt^{\lie g}[0]}^{\otimes k}\cong \bt^{\lie g}[k],\qquad (\bt_\Psi^{\lie g}(F)[1])_{\bt^{\lie g}[0]}^{\otimes k}\cong \bt_\Psi^{\lie g}(F)[k],\qquad k\ge 1.
$$
The first statement of the following proposition is now immediate while the second follows by using ~\lemref{SALG20}.

\begin{prop}The map $\bom: T_{\bt^{\lie g}[0]}(\bt^{\lie g}[1])\to \bt^{\lie g}$ is an isomorphism of
 $\bz_+$-graded associative algebras. In particular,
the algebras $\bs^{\lie g}$, $\be^{\lie g}$ are quadratic. Similarly, for $F\subset P^+$ interval closed the
restriction $\bom: T_{\bt^{\lie g}_\Psi(F)[0]}(\bt^{\lie g}_\Psi(F)[1])\to \bt_\Psi^{\lie g}(F)$
is an isomorphism of $\bz_+$-graded associative algebras. In particular,
the algebras $\bs_\Psi^{\lie g}(F)$, $\be^{\lie g}_\Psi(F)$ are quadratic.
\hfill\qedsymbol
\end{prop}

\subsection{}\label{SALG90}

\begin{prop}
Let $\lambda\in P^+$ and assume that~$\Psi=\Psi(\psi)$ for some~$\psi\in P$.
The algebras~$\bs_\Psi^{\lie g}$ and~$\bs_\Psi^{\lie g}(\lambda\le_\Psi)$ are Koszul.

\end{prop}
\begin{pf} We use the notation of Section~\ref{SALG0}. Let $A=\bs_\Psi^{\lie g}$ and for $F\subset P^+$ set $A(F)=\bs_\Psi^{\lie g}(F)$. Clearly $A=A(P^+)$.  In~\propref{CAT110} we proved that the algebra $A(\le_\Psi\mu)$ is Koszul for all $\mu\in P^+$. Hence the Koszul complex of this algebra is exact. We now describe the relationship between the Koszul complex of $A$ and that of $A(\le_\Psi\mu)$, $\mu\in P^+$.

Setting  $N^r=N^r_A$ and $N^r(F)=N^r_{A(F)}$, we see that\begin{gather*}  N^r 1_\mu=\bigoplus_{\nu\le_\Psi\mu} 1_\nu N^r 1_\mu\ \cong\  N^r(\le_\Psi\mu),\\ N^r=\bigoplus_{\mu\in P^+} N^r 1_\mu\ \cong \bigoplus_{\mu\in P^+}N^r(\le_\Psi\mu).\end{gather*}

This gives,$$ A\tensor_{A[0]} N^r 1_\mu=\bigoplus_{\nu\le_\Psi\mu} A 1_\nu\tensor_{A[0]}
1_\nu N^r 1_\mu=\bigoplus_{\nu,\xi\le_\Psi\mu}A 1_\nu\tensor_{A[0]}
1_\xi N^r 1_\mu,$$
which in turn implies that
$$A\tensor_{A[0]} N^r 1_\mu\cong A(\le_\Psi\mu)\tensor_{A(\le_\Psi\mu)[0]} N^r(\le_\Psi\mu).$$ Moreover this isomorphism is compatible with the maps $\partial_r$ and hence we get
$$
\ker \partial_r|_{A\tensor_{A[0]} N^r 1_\mu}\subset \partial_{r+1}(A\tensor_{A[0]} N^{r+1} 1_\mu),
$$
which proves that the Koszul complex for $A$ is exact. The proof when $F$ is the set $\lambda\le_\Psi$ is similar and we omit the details.
\end{pf}

\section{Quadratic dual of the algebra~\texorpdfstring{$\bs_\Psi^{\lie g}$}{\bfseries S\_Psi g}}
There are two notions of quadratic duals which appear in the literature. One definition can be found in the study of Koszul quotients of path algebras of quivers. The other one given in \cite[Definition~2.8.1]{BGS} applies to quadratic algebras $A$ which satisfy the additional condition that $A[r]$, for all $r\ge 0$ is a finitely generated left $A[0]$-module. In our case   the two definitions are equivalent for the algebras $\bs_\Psi^{\lie g}(F)$ when $F$ is finite and interval closed but only the first definition can be used to define the quadratic dual of $\bs_\Psi^{\lie g}$.

\subsection{}\label{QUAD30}
Fix~$\Psi=\Psi(\psi)$, $\psi\in P$.
\begin{lem}
Let~$F\subset P^+$ be interval closed. Let~$\ba$ be $\bt$, $\bs$ or~$\be$. Then $(\ba_\Psi^{\lie g}(F))^{op}$ is isomorphic to the
the subalgebra
$$
\bigoplus_{\lambda\le_\Psi\mu\in F} 1_\mu \ba^{\lie g}[d_\Psi(\lambda,\mu)] 1_\lambda\subset \ba^{\lie g},
$$ of $\ba_\Psi$.
\end{lem}
\begin{pf} Note that $A[k]^*\cong A[k]$ as $\lie g$-modules if  $A$ is one of  $T(\lie g)$, $S(\lie g)$ or $\bigwedge\lie g$. Moreover for $M,N\in\Ob\cal F(\lie g)$, we have $(M^*\otimes N)^{\lie g}\cong(N^*\otimes M)^{\lie g}$. Furthermore, if~$K\in\Ob\cal F(\lie g)$, then
the canonical map (cf.~\lemref{PRE10})
$$
(M^*\tensor N)^{\lie g}\tensor (N^*\tensor K)^{\lie g}\to (M^*\tensor K)^{\lie g}
$$
induces the canonical map
$$
(K^*\tensor N)^{\lie g}\tensor (N^*\tensor M)^{\lie g}\to (K^*\tensor M)^{\lie g}.
$$
Let~$\lambda\le_\Psi\mu\in F$. By the above we
have an isomorphism of vector spaces
\begin{align*} 1_\lambda\ba^{\lie g}[d_\Psi(\lambda,\mu)]1_\mu&
=(V(\lambda)^*\tensor A[d_\Psi(\lambda,\mu)]\tensor V(\mu))^{\lie g}\\
&\cong (V(\mu)^*\tensor A[d_\Psi(\lambda,\mu)]\tensor V(\lambda))^{\lie g}
=1_\mu\ba^{\lie g}[d_\Psi(\lambda,\mu)]1_\lambda
\end{align*}
which extends to the isomorphism of algebras
\begin{equation*}
\bigoplus_{\lambda\le_\Psi\mu\in F} 1_\mu \ba^{\lie g}[d_\Psi(\lambda,\mu)] 1_\lambda\cong (\ba^{\lie g}_\Psi(F))^{op}.
\qedhere
\end{equation*}
\end{pf}

\subsection{}\label{QUAD5}

Let $(\cdot,\cdot)_{\lie g}$ be the Killing form of~$\lie g$. Define
$
(\cdot,\cdot)_{T(\lie g)}: T(\lie g)\tensor T(\lie g)\to \bc
$
by extending linearly the assignment \begin{gather*} (T^r(\lie g),T^s(\lie g))=0,\qquad r\not=s,\\
(x_1\tensor \cdots \tensor x_r,y_r\tensor\cdots\tensor  y_1)_{T(\lie g)}=(x_1,y_1)_{\lie g} \cdots (x_r,y_r)_{\lie g},
\qquad x_i,y_i\in\lie g,\, 1\le i\le r.
\end{gather*}
Then $(\cdot,\cdot)_{T(\lie g)}$ is a $\lie g$-invariant  symmetric bilinear form on $T(\lie g)$.
It is easy to check that the restrictions of~$(\cdot,\cdot)_{T(\lie g)}$ to
$S^2(\lie g)\tensor S^2(\lie g)$ and $\bigwedge^2 \lie g\tensor \bigwedge^2 \lie g$ are non-degenerate.
Since $T^2(\lie g)=S^2(\lie g)\oplus \bigwedge^2 \lie g$ and
$$
(x,y)_{T(\lie g)}=0,\qquad x\in S^2(\lie g),\, y\in \bigwedge^2 \lie g,
$$
it follows that
\begin{equation}\label{QUAD5.10}
\{ x\in T^2(\lie g)\,:\, (x,\bigwedge^2 \lie g)_{T(\lie g)}=0\}=S^2(\lie g).
\end{equation}

\subsection{}\label{QUAD10}
Define
$$
\langle\cdot,\cdot\rangle: \bt\tensor \bt\to \bc
$$
 by
$$
\langle f\tensor a\tensor v,f'\tensor b\tensor v'\rangle=(a,b)_{T(\lie g)} f'(v) f(v'),\qquad v,v'\in\bv,\,f,f'\in\bv^{\circledast}, a,b\in T(\lie g).
$$
It is easy to check that $\langle\cdot,\cdot\rangle$ is a symmetric non-degenerate and  $\lie g$-invariant form. Moreover,
\begin{gather} \label{inv}\langle 1_\lambda ua,v 1_\mu\rangle=\delta_{\lambda,\mu}\langle u,av\rangle,\qquad u,v\in\bt,\, a\in\bt[0],\, \lambda,\mu\in P^+.
\end{gather}
\begin{lem}
The  restriction of $\langle\cdot,\cdot\rangle$  to~$\bt^{\lie g}\tensor \bt^{\lie g}$ is non-degenerate. In particular, the restriction of $\langle\cdot,\cdot\rangle$ to $\bt^{\lie g}_\Psi(F)\tensor (\bt^{\lie g}_\Psi(F))^{op}$ is non-degenerate.
\end{lem}

\begin{pf}
Since $$\bt^{\lie g}=\bigoplus_{\lambda,\mu\in P^+,k\in\bz_+}
1_\lambda \bt^{\lie g}[k] 1_\mu
$$
and $\langle \bt[k],\bt[s]\rangle=0$, $k\ne s$, it suffices to show that for all $\lambda,\mu\in P^+$,
\begin{equation}\label{nondeg} x\in  1_\mu\bt^{\lie g}[k]1_\lambda,\qquad
 \langle x,1_\mu\bt^{\lie g}[k]1_\lambda\rangle=0
\implies \langle x,\bt[k]\rangle =0.
\end{equation}

Let $y\in 1_\nu\bt[k]1_\xi$.  If $\nu\ne\mu$ or $\xi\ne\lambda$ then \eqref{nondeg} follows from~\eqref{inv}. Otherwise, write $y=\sum_{\zeta\in P^+}y_\zeta$, where $y_\zeta$ is in the $\zeta$-isotypical component of
$ 1_\mu\bt[k]1_\lambda\cong V(\mu)^*\otimes T(\lie g)\otimes V(\lambda)$. Since the form is $\lie g$-invariant, it follows now that $\langle x,y_\zeta\rangle=0$ if $\zeta\ne 0$. This proves \eqref{nondeg}. The second assertion is immediate by~\lemref{QUAD30}.
\end{pf}

\subsection{}\label{QUAD50}
{} Let~$F\subset P^+$ be interval closed.
For~$\ba=\bs$ or~$\be$,
let $$\br_{\ba}=\ker(\bt^{\lie g}_\Psi(F)\to \ba^{\lie g}_\Psi(F))\cap \bt^{\lie g}[2],\ \ \br_\ba^!=\{ x\in \bt^{\lie g}[2]\,:\, \langle \br_\ba,x\rangle=0\},$$ and set
$$(\ba^{\lie g}_\Psi(F))^!=(\bt^{\lie g}_\Psi(F))^{op}/\langle\br_{\ba}^!\rangle.
$$
If~$F$ is finite, it is not hard to see that this algebra is isomorphic to the right quadratic dual of~$\ba_\Psi^{\lie g}(F)$ as defined in~\cite{BGS}.
\begin{prop}
Let~$F\subset P^+$ be interval closed. Then $(\bs^{\lie g}_\Psi(F))^!\cong (\be^{\lie g}_\Psi(F))^{op}$.
In particular, for all~$\lambda'\le_\Psi\lambda\in P^+$, the algebras $\be^{\lie g}_\Psi(\le_\Psi\lambda)$, $\be^{\lie g}_\Psi([\lambda',\lambda]_\Psi)$
are Koszul.
\end{prop}
\begin{pf}
To prove the first assertion, recall that by~\lemref{SALG20}\eqref{SALG20.ii}
$$
\br_{\bs}=\bigoplus_{\lambda\le_\Psi\mu\in F\,:\,d_\Psi(\lambda,\mu)=2} (V(\lambda)^*\tensor \bigwedge^2\lie g\tensor V(\mu))^{\lie g}.
$$
Using~\lemref{QUAD10} and~\eqref{QUAD5.10},
we conclude that
$$
\br_{\bs}^!=\bigoplus_{\lambda\le_\Psi\mu\in F\,:\,d_\Psi(\lambda,\mu)=2} (V(\mu)^*\tensor S^2(\lie g)\tensor V(\lambda))^{\lie g}
$$
It remains to apply~\lemref{SALG20}\eqref{SALG20.ii} and~\lemref{QUAD30}.
The second assertion follows immediately from~\cite[Proposition~2.2.1]{BGS}.
\end{pf}

\subsection{}\label{QUAD70}
The following proposition completes the proof of~\thmref{kosalg}.
\begin{prop}
Let~$\lambda\in P^+$.
The algebras~$\be^{\lie g}_\Psi$, $\be^{\lie g}_\Psi(\lambda\le_\Psi)$ are Koszul.
\end{prop}
\begin{pf}
As in~\propref{SALG90}, it suffices to prove that the Koszul complex for~$\be^{\lie g}_\Psi$ is exact. The argument is similar to
that in~\propref{SALG90}, with the algebras~$\bs$ replaced by the corresponding algebras~$\be$ and is omitted.
\end{pf}

\end{document}